\documentstyle{amsppt}

\NoBlackBoxes

\define\Lip {\operatorname{Lip}}
\define\ep {\varepsilon}
\define\End {\operatorname{End}}
\define\Dbar {\overline{D}}
\define\Ubar {\overline{U}}
\define\Obar {\overline{\Cal{O}}}
\define\Sbar {\overline{\Cal{S}}}
\define\zbar {\overline{z}}
\define\zetabar {\overline{\zeta}}
\define\wtB {\widetilde{B}}
\redefine\Re {\operatorname{Re}}
\redefine\Im {\operatorname{Im}}

\define\RR {\Bbb R}
\define\CC {\Bbb C}
\define\ZZ {\Bbb Z}
\define\pa {\partial}

\magnification=1200
\vsize=8.5truein

\centerline{\bf The Complex Frobenius Theorem for Rough Involutive Structures}
\vskip 15pt
\centerline{\smc {C.~Denson Hill} and {Michael Taylor}
\footnote{{\it 2000 Mathematics Subject Calssification.} Primary 35N10 
\newline 
The second author was partially supported by NSF grant DMS-0139726}
}

$$\text{}$$
{\smc Abstract.}  We establish a version of the complex Frobenius theorem in 
the context of a complex subbundle 
$\Cal{S}$ of the complexified tangent bundle
of a manifold, having minimal regularity.  If the subbundle $\Cal{S}$ defines
the structure of a Levi-flat CR-manifold, it suffices that $\Cal{S}$ be
Lipschitz for our results to apply.  A principal tool in the analysis is a
precise version of the Newlander-Nirenberg theorem with parameters, for 
integrable almost complex structures with minimal regularity, which builds on
previous recent work of the authors.

$$\text{}$$
{\smc Contents}:
1. Introduction,
2. Real Frobenius theorem for involutive Lipschitz bundles,
3. The pull-back of a Levi-flat CR structure,
4. The Newlander-Nirenberg theorem with parameters,
5. Structure of Levi-flat CR-manifolds,
6. The complex Frobenius theorem,
A. A Frobenius theorem for real analytic, complex vector fields,
B. The case of two-dimensional leaves.

$$\text{}$$
{\bf 1. Introduction}
\newline {}\newline

The complex Frobenius theorem elucidates the structure of a complex 
subbundle $\Cal{S}$ of the complexified tangent bundle $\CC T\Omega$
of a smooth manifold $\Omega$, satisfying an involutivity condition, 
which can be stated as follows: if $X$ and $Y$ are (sufficiently regular)
sections of $\Cal{S}$, then
$$
[X,Y]\ \text{ is a section of }\ \Cal{S},
\tag{1.1}
$$
and
$$
[X,\overline{Y}]\ \text{ is a section of }\ \Cal{S}+\Sbar.
\tag{1.2}
$$
Here, as usual, if $X=X_0+iX_1$ and $X_0,X_1$ are real vector fields, we
write $\overline{X}=X_0-iX_1$, and the fiber of $\Sbar$ over $p\in\Omega$
is given as
$$
\Sbar_p=\{u-iv:u+iv\in \Cal{S}_p,\ u,v\in T_p\Omega\}.
\tag{1.3}
$$
We also assume $\Cal{S}+\Sbar$ is a subbundle of $\CC T\Omega$.

In case $\Cal{S}=\CC \Cal{S}_0$ is the complexification of a subbundle 
$\Cal{S}_0\subset T\Omega$, the condition (1.1) just says $\Cal{S}_0$ is
involutive. (Here, $\Cal{S}=\Sbar$, and (1.2) provides no additional 
constraint.)  In this case the result reduces to the real Frobenius theorem.

An opposite extreme arises when $\Omega$ has an almost complex structure,
a section $J$ of $\End T\Omega$ satisfying $J^2=-I$ (which implies that 
dim $\Omega$ is even).  We set
$$
\Cal{S}_p=\{u+iJu:u\in T_p\Omega\},
\tag{1.4}
$$
so a section of $\Cal{S}$ has the form $X+iJX$, for a general real vector
field $X$.  The condition (1.1) is that if also $Y$ is a real vector field,
then $[X+iJX,Y+iJY]=Z+iJZ$ for a real vector field $Z$.  This is equivalent
to the vanishing of the Nijenhuis tensor, defined by
$$
N(X,Y)=[X,Y]-[JX,JY]+J[X,JY]+J[JX,Y].
\tag{1.5}
$$
The content of the Newlander-Nirenberg theorem [NN] is that under this 
formal integrability hypothesis $\Omega$ has local holomorphic coordinates,
i.e., functions $u_1,\dots,u_k:\Cal{O}\rightarrow\CC$ forming a coordinate
system on a neighborhood $\Cal{O}$ of a given $p\in\Omega$, such that 
$(X+iJX)u_\ell\equiv 0$ for all real vector fields $X$.  Thus $\Omega$ has
the structure of a complex manifold.  In this case, $\Cal{S}+\Sbar=\CC T
\Omega$, so (1.2) automatically holds.  There are other cases, where (1.2)
has a nontrivial effect, as will be seen below.  

The complex Frobenius theorem was established in [Ni] for $C^\infty$ 
bundles $\Cal{S}\subset\CC T\Omega$ satisfying (1.1)--(1.2).  A major
ingredient in the proof was the Newlander-Nirenberg theorem, which had
been established in [NN] for almost complex structures with a fairly high
degree of smoothness.  Later proofs of the Newlander-Nirenberg theorem,
by [NW] and by [M], work for almost complex structures $J$ of class $C^{1+r}$
with $r>0$, i.e., when $J$ has H{\"o}lder continuous first order derivatives.
In [HT] the needed regularity on $J$ was reduced to $J\in C^r$ with $r>1/2$.
(More general conditions were considered in [HT], which we will not discuss 
here.)  The case of Lipschitz $J$ found an immediate application in [LM].

Regarding the real Frobenius theorem, standard arguments, though frequently
phrased in the context of smooth subbundles of $T\Omega$, work for $C^1$
bundles.  The real Frobenius theorem was extended in [Ha] to include
Lipschitz subbundles.

Our main goal here is to extend Nirenberg's complex Frobenius theorem to the
setting of rough bundles $\Cal{S}\subset \CC T\Omega$ satisfying (1.1)--(1.2).
We will assume that $\Cal{S}$ and $\Cal{S}+\Sbar$ are Lipschitz 
subbundles of $\CC T\Omega$.  Note that if $X$ and $Y$ are Lipschitz
sections of $\Cal{S}$, then $[X,Y]$ and $[X,\overline{Y}]$ are vector
fields with $L^\infty$ coefficients.  For an important class of bundles 
$\Cal{S}$, namely those giving rise to Levi-flat CR-structures (defined below)
this regularity hypothesis will suffice.  In the general case we need an 
additional hypothesis, given in (1.16) below.  We mention that [Ho] 
established a version of a complex Frobenius theorem in a setting of $C^1$ 
vector fields, with $C^1$ commutators, but with a somewhat different thrust.

We now set up a basic strategy for obtaining such a complex Frobenius 
theorem, and indicate what extra analysis has to be done to treat the 
non-smooth case.  It is convenient to begin by constructing some further 
subbundles of the real tangent bundle $T\Omega$.  For each $p\in\Omega$, set
$$
\aligned
\Cal{E}_p&=\{u\in T_p\Omega:u+iv\in\Cal{S}_p,\text{for some }v\in T_p\Omega\} 
\\ &=\{w+\overline{w}:w\in\Cal{S}_p\},
\endaligned
\tag{1.6}
$$
the fiber over $p$ of a Lipschitz bundle $\Cal{E}$.  Noting that if $u,v
\in T_p\Omega$ and $u+iv\in\Cal{S}_p$, then also $v-iu\in \Cal{S}_p$, so
$v\in\Cal{E}_p$, we see that
$$
\Cal{S}+\Sbar=\CC \Cal{E}.
\tag{1.7}
$$
Next, set
$$
\Cal{V}_p=\Cal{S}_p\cap T_p\Omega,
\tag{1.8}
$$
the fiber over $p$ of a Lipschitz vector bundle $\Cal{V}$.  Note that if 
$u,v\in T_p\Omega$,
$$
\aligned
u+iv\in \Cal{S}_p\cap\Sbar_p&\Longleftrightarrow u+iv\in\Cal{S}_p\ 
\text{ and }\ u-iv\in\Cal{S}_p \\
&\Longleftrightarrow u\in\Cal{S}_p\ \text{ and }\ v\in\Cal{S}_p.
\endaligned
\tag{1.9}
$$
Hence
$$
\Cal{S}\cap\Sbar=\CC \Cal{V}.
\tag{1.10}
$$
The hypotheses (1.1)--(1.2) imply $\Cal{E}$ and $\Cal{V}$ are involutive
subbundles of $T\Omega$, i.e.,
$$
\aligned
X,Y\in\Lip(\Omega,\Cal{E})&\Longrightarrow [X,Y]\in L^\infty(\Omega,\Cal{E}), 
\\ X,Y\in\Lip(\Omega,\Cal{V})&\Longrightarrow [X,Y]\in
L^\infty(\Omega,\Cal{V}).
\endaligned
\tag{1.11}
$$
On the other hand, one does not recover (1.1)--(1.2) from (1.11) alone, as
our second example illustrates.  In that example, with $\Cal{S}_p$ given by
(1.4), we have $\Cal{E}=T\Omega,\ \Cal{V}=0$, and (1.11) always holds, 
regardless of whether $N$ in (1.5) vanishes.  To capture (1.1)--(1.2),
an additional structure arises.

Namely, one has a complex structure on the quotient 
bundle $\Cal{E}/\Cal{V}$, defined as follows.  Take $u\in\Cal{E}_p$, 
so there exists $v\in T_p\Omega$ such that $u+iv\in\Cal{S}_p$; in fact,
$v\in\Cal{E}_p$.  We propose to set $Ju=v$, so the element of $\Cal{S}_p$
has the form $u+iJu$.  However, the element $v$ associated to $u\in 
\Cal{E}_p$ is not necessarily unique.  In fact, given $u,v,v'\in T_p\Omega$
and $u+iv\in \Cal{S}_p$, we have
$$
u+iv'\in\Cal{S}_p\Leftrightarrow i(v-v')\in\Cal{S}_p\Leftrightarrow
v-v'\in\Cal{S}_p\cap\Cal{E}_p=\Cal{V}_p.
\tag{1.12}
$$
In other words, given $u\in\Cal{E}_p$, the residue class of $Ju$ is well 
defined in $\Cal{E}_p/\Cal{V}_p$.  Furthermore, if $u\in\Cal{V}_p$, one can
take $v=0$, so $J$ descends from a linear map $\Cal{E}_p\rightarrow 
\Cal{E}_p/\Cal{V}_p$ to
$$
J_p:\Cal{E}_p/\Cal{V}_p\longrightarrow \Cal{E}_p/\Cal{V}_p,
\tag{1.13}
$$
yielding
$$
J\in\Lip(\Omega,\End\, \Cal{E}/\Cal{V}).
\tag{1.14}
$$
Since $u+iv\in\Cal{S}_p\Leftrightarrow v-iu\in\Cal{S}_p$, we also have $J^2=
-I$.  The integrability hypotheses (1.1)--(1.2) are equivalent to (1.11), 
coupled to an integrability hypothesis on $J$, which we describe below.

Let us first consider the case $\Cal{V}=0$.  Then $J$ is a complex structure 
on the involutive bundle $\Cal{E}$, and (generalizing (1.4)) we have
$$
\Cal{S}_p=\{u+iJu:u\in\Cal{E}_p\},
\tag{1.15}
$$
or equivalently Lipschitz sections of $\Cal{S}$ have the form $X+iJX$, where
$X$ is a Lipschitz section of $\Cal{E}$.  Then the involutivity hypothesis 
(1.1)--(1.2) is equivalent to the involutivity of $\Cal{E}$ plus the 
vanishing of $N$, given by (1.5), for $X,Y\in\Lip(\Omega,\Cal{E})$.  
One says that $\Omega$ has the structure of a Levi-flat CR manifold.  The 
real Frobenius theorem implies that $\Omega$ is foliated by leaves tangent to 
$\Cal{E}$.  Each such leaf then inherits an almost complex structure, and
the Newlander-Nirenberg theorem implies each such leaf has local holomorphic
coordinates.  Briefly put, $\Omega$ is foliated by complex manifolds.  The
complex Frobenius theorem in this context says a little more.  Namely, any
$p\in\Omega$ has a neighborhood $\Cal{O}$ on which there are functions
$u_1,\dots,u_k$, providing holomorphic coordinates on each leaf, intersected 
with $\Cal{O}$, and having some regularity on $\Cal{O}$.  In the case of a
$C^\infty$ bundle $\Cal{S}$, [Ni] obtained such $u_j\in C^\infty(\Cal{O})$.
In the context of Lipschitz structures, we obtain certain H{\"o}lder 
continuity of $u_j$, described in further detail below.  A key ingredient in 
the analysis is a Newlander-Nirenberg theorem with parameters.  In the smooth
case this follows by the methods of [NN], as noted there and used in [Ni].
We devote \S{4} to a consideration of families of integrable almost complex 
structures with minimal regularity, building on techniques of [M] and of 
[HT].

We now turn to the case $\Cal{V}\neq 0$.  In this case, we supplement the 
Lipschitz hypotheses on $\Cal{S}$ and $\Cal{S}+\Sbar$ with the following
hypothesis.  Say $\text{dim}\, \Cal{V}_p=\ell\le k=\text{dim}\, \Cal{E}_p$.
We assume that each $p\in\Omega$ has a neighborhood on which there is a 
local Lipschitz frame field $\{X_1,\dots,X_k\}$ for $\Cal{E}$, such that
$\{X_1,\dots,X_\ell\}$ is a local frame field for $\Cal{V}$ and
$$
[X_i,X_j]=0,\quad 1\le i,j\le k.
\tag{1.16}
$$
This can be regarded as an hypothesis on the regularity with which $\Cal{V}$
sits in $\Cal{E}$; we discuss it further in \S{6}.  We will show that
$$
J\ \text{ is invariant under }\ \Cal{F}^t_{X_i},\ \ 1\le i\le \ell,
\tag{1.17}
$$
where $\Cal{F}^t_{X_i}$ is the flow generated by $X_i$.  Hence
we can mod out by the $\Cal{F}^{t_1}_{X_1}\circ\cdots\circ
\Cal{F}^{t_\ell}_{X_\ell}$ action, to obtain
$$
\pi:\Omega\longrightarrow M,
\tag{1.18}
$$
(perhaps after localizing), and on $M$ we have a Levi-flat CR structure.
Leafwise holomorphic functions on (open subsets of) $M$ pull back to 
functions on (open subsets of) $\Omega$, and results on their existence and 
regularity essentially constitute the complex Frobenius theorem for the 
bundle $\Cal{S}$.

The rest of this paper is organized as follows.  Section 2 treats the real
Frobenius theorem for involutive Lipschitz bundles.  We recall some results
of [Ha] and establish some further results, regarding the regularity of the
diffeomorphism constructed to flatten out the leaves of the foliation.
In \S{3} we consider Levi-flat CR manifolds, in the Lipschitz
category, even allowing for rougher $J$, and examine how such a structure 
pulls back under a leaf-flattening diffeomorphism from \S{2}, to yield a
parametrized family of manifolds carrying integrable almost complex 
structures.  This sets us up for a study of the Newlander-Nirenberg theorem 
with parameters, which we carry out in \S{4}.  

In \S{5} we tie together the material of \S\S{2--4} to obtain results on 
the existence and regularity of
functions on open sets of a Lipschitz Levi-flat CR manifold $\Omega$ that are
leafwise holomorphic (functions known as CR functions).  Our primary result,
Proposition 5.1, yields CR functions $\varphi_j,\ 1\le j\le m+n-k$, on a 
neighborhood $U_1$ of a point $p\in\Omega$, having the property that
$$
\Phi=(\varphi_1,\dots,\varphi_{m+n-k}):U_1\longrightarrow \CC^m\times
\RR^{n-k}
\tag{1.19}
$$
is a homeomorphism of $U_1$ onto an open subset, and such that, given $s<1/2$,
$\varphi_j$ and $X\varphi_j$ are H{\"o}lder continuous of degree $s$, for any
$X\in\Lip(U_1,\Cal{E})$.  A complementary result, Proposition 5.2, shows that 
$\Phi$ in (1.19) can be taken to be
a $C^1$ diffeomorphism, provided that $\Cal{S}$, and hence
$\Cal{E}$ and $J$, are regular of class $C^\rho$ for some $\rho>3/2$.  
The results of [HT] extending the Newlander-Nirenberg theorem to cases 
where the almost complex structure is merely $C^{1/2+\ep}$ regular, and
complementary results of \S{4}, play an important role in the proof.
We end \S{5} with a brief discussion of $C^{1,1}$ submanifolds of $\CC^N$
that have the structure of Levi-flat CR-manifolds.
The general complex Frobenius theorem is then treated in \S{6}.

At the end of this paper we have two appendices.  Appendix A is devoted to
a Frobenius theorem for real analytic, complex vector fields.  There are 
classical results of this nature; cf.~[Ni] for some references.  
One motivation for us to include a self contained 
treatment of such a result here arises 
from the nature of our analysis of the Newlander-Nirenberg theorem with
parameters in \S{4}.  Following [M], we construct the local holomorphic
coordinate chart as a composition, $F=G\circ H$.  The map $H$ is obtained via
an implicit function theorem, the use of which enables us to keep track of its
dependence on a parametrized family of integrable almost complex structures.
The construction of $H$ arranges things so that constructing $G$ amounts to
establishing the Newlander-Nirenberg theorem in the real analytic category,
a task to which the material of Appendix A is applicable, and this material 
makes it clear how the factor $G$ depends on parameters.

Finally, Appendix B gives a special treatment of the construction of CR 
functions on a rough Levi-flat CR manifold whose leaves have real dimension
2.  The classical method of constructing isothermal coordinates is adapted
to this problem and yields sharper results than one obtains in the case 
of higher dimensional leaves via the methods of \S{4}.  This leads to 
improved results in \S{5} in the case of 2-dimensional leaves, as is noted 
there.

We end this introduction with a few remarks on function spaces arising in
our analysis.  For a smoothly bounded domain $\Ubar,\ C^r(\Ubar)$ denotes the
space of functions with derivatives of order $\le r$ continuous on $\Ubar$,
if $r$ is a positive integer.  If $r=k+s,\ k\in\ZZ^+,\ 0<s<1$, it denotes the 
space of functions whose $k$th order derivatives are H{\"o}lder continuous 
of order $s$.  In addition, we make use of Zygmund spaces $C^r_*(\Ubar)$, 
coinciding with $C^r(\Ubar)$ for $r\in\RR^+\setminus\ZZ^+$, and having nice
interpolation properties at $r\in\ZZ^+$.  The spaces $C^r_*(\Ubar)$ are also
defined for $r<0$.  There are a number of available treatments of Zygmund 
spaces; we mention Chapter 13, \S{8} of [T] as one source.  As is usual, 
$\Lip(\Ubar)$ denotes the space of Lipschitz continuous functions, i.e., 
functions H{\"o}lder continuous of exponent one, and $C^{1,1}(\Ubar)$ denotes
the space of functions whose first order derivatives belong to $\Lip(\Ubar)$.

$$\text{}$$
{\bf 2. Real Frobenius theorem for involutive Lipschitz bundles}
\newline {}\newline

Let $\Cal{E}$ be a sub-bundle of the tangent bundle $T\Omega$, of fiber
dimension $k$.  We assume $\Cal{E}$ is Lipschitz, in the sense that
any $p_0\in\Omega$ has a neighborhood $\Cal{O}$ on which there are Lipschitz
vector fields $X_1,\dots,X_k$ spanning $\Cal{E}$ at each point.  We make the
involutivity hypothesis that $[X_i,X_j]$ is a section of $\Cal{E}$ at 
almost all points of $\Cal{O}$, or equivalently that there exist
$c^\ell_{ij}\in L^\infty(\Cal{O})$ such that
$$
[X_i,X_j]=\sum\limits_\ell c^\ell_{ij}(x)\, X_\ell.
\tag{2.1}
$$
We want to discuss the existence and qualitative properties of the foliation
of $\Omega$ whose leaves are tangent to $\Cal{E}$.

We may as well assume $k<n=\text{dim}\, \Omega$.  Suppose we have coordinates
centered at $p_0$ such that $X_j(p_0)$ form the first $k$ standard basis
elements of $\RR^n$, for $1\le j\le k$.  If we denote by $\widetilde{X}_j(x)$
the image of $X_j(x)$ under the standard projection $\RR^n\rightarrow \RR^k$,
we have
$$
\widetilde{X}_i(x)=\sum\limits_{j=1}^k A_{ij}(x)\, \pa_j,
\tag{2.2}
$$
with $A_{ij}\in\Lip(\Cal{O}),\ A_{ij}(p_0)=\delta_{ij}$, hence $(A_{ij}(x))$
an invertible $k\times k$ matrix, with inverse $(B_{ij}(x))$, for $x$ in a
neighborhood of $p_0$ (which we now denote $\Cal{O}$).  We set
$$
Y_i=\sum\limits_j B_{ij}(x)\, X_j,\quad 1\le i\le k.
\tag{2.3}
$$
It follows that
$$
Y_i=\pa_i+Y^\#_i,\quad Y^\#_i=\sum\limits_{\ell\ge k+1} D_{i\ell}(x)\pa_\ell,
\quad 1\le i\le k.
\tag{2.4}
$$
Also (2.1) implies
$$
[Y_i,Y_j]=\sum\limits_\ell \tilde{c}^\ell_{ij}(x)\, Y_\ell,
\tag{2.5}
$$
for certain $\widetilde{c}^\ell_{ij}\in L^\infty(\Cal{O})$.  Comparison of 
(2.4) and (2.5) yields $\tilde{c}^\ell_{ij}\equiv 0$, so we have a local
Lipschitz frame field for $\Cal{E}$ satisfying
$$
[Y_i,Y_j]=0,\quad 1\le i,j\le k.
\tag{2.6}
$$
The key result on the existence of a foliation tangent to $\Cal{E}$ is the 
following result of [Ha].

\proclaim{Proposition 2.1} Let $Y_j$ be Lipschitz vector fields on $\Cal{O}$
satisfying (2.6).  For any compact $K\subset\Cal{O}$ there exists
$\delta>0$ such that there is a unique solution $y=y(t,x_0)=
y(t_1,\dots,t_k,x_0)$ to
$$
\frac{\pa y}{\pa t_j}=Y_j(y),\quad 1\le j\le k,\quad y(0,x_0)=x_0,
\tag{2.7}
$$
given $x_0\in K,\ |t_j|<\delta$.  Furthermore, $y(t,x)$ is Lipschitz in
$(t,x)$.
\endproclaim

In fact, this result is a special case of Corollary 4.1 of [Ha].  We make
some further comments on it.  If $\Cal{F}^t_{Y_j}$ denotes the flow
generated by $Y_j$, we see that
$$
y(0,\dots,0,t_k,x)=\Cal{F}^{t_k}_{Y_k}(x).
\tag{2.8}
$$
Then
$$
y(0,\dots,0,t_{k-1},t_k,x)=\Cal{F}^{t_{k-1}}_{Y_{k-1}}\circ
\Cal{F}^{t_k}_{Y_k}(x),
\tag{2.9}
$$
and inductively
$$
y(t_1,\dots,t_k,x)=\Cal{F}^{t_1}_{Y_1}\circ\cdots\circ\Cal{F}^{t_k}_{Y_k}(x).
\tag{2.10}
$$
The order can be changed, and we have
$$
\Cal{F}^{t_i}_{Y_i}\circ \Cal{F}^{t_j}_{Y_j}(x)=
\Cal{F}^{t_j}_{Y_j}\circ\Cal{F}^{t_i}_{Y_i}(x),
\tag{2.11}
$$
for $x\in K,\ |t_i|,|t_j|<\delta$.  

Conversely, once one knows that (2.11) follows from (2.6), one can prove 
Proposition 2.1.  However, this implication is less straightforward for 
Lipschitz vector fields than it is for smooth vector fields.  In connection
with this, we mention the following analytical point, which plays a key
role in the proof in [Ha].  Namely, let $\{J_\ep:0<\ep\le 1\}$ be a
Friedrichs mollifier and let $Y_i,Y_j$ be Lipschitz vector fields satisfying
(2.6).  Then, as $\ep\rightarrow 0$,
$$
[J_\ep Y_i,J_\ep Y_j]\longrightarrow 0,
\tag{2.12}
$$
locally uniformly on $\Cal{O}$.  Actually this is a reformulation (of a 
special case) of Proposition 5.3 of [Ha].  It is stronger and more useful
than the obvious fact that such convergence holds $\text{weak}^*$ in 
$L^\infty$.  What is behind it is the more general fact that, for any two
Lipschitz vector fields $X$ and $Y$ on $\Cal{O}$, 
$$
[J_\ep X,J_\ep Y]-J_\ep [X,Y]\longrightarrow 0,
\tag{2.13}
$$
locally uniformly on $\Cal{O}$.  This follows from the fact that
$$
f\in\Lip(\Cal{O}),\ g\in L^\infty(\Cal{O})\Longrightarrow
(J_\ep f)(J_\ep g)-J_\ep(fg)\rightarrow 0,
\tag{2.14}
$$
locally uniformly on $\Cal{O}$, and since clearly $J_\ep f\rightarrow f$
locally uniformly on $\Cal{O}$ this in turn is equivalent to the fact that
$$
f\in\Lip(\Cal{O}),\ g\in L^\infty(\Cal{O})\Longrightarrow
f\, J_\ep g-J_\ep(fg)\rightarrow 0,
\tag{2.15}
$$
locally uniformly on $\Cal{O}$, which is a standard Friedrichs-type 
commutator estimate.

We record that $y(t,x)$ has extra regularity in $t$.  
\proclaim{Corollary 2.2} For each $j\in\{1,\dots,k\}$,
$$
\frac{\pa}{\pa t_j} y(t,x)\ \text{ is Lipschitz in }\ (t,x).
\tag{2.16}
$$
\endproclaim
\demo{Proof} Clearly the right side of (2.7) is Lipschitz in $(t,x)$.
\enddemo

Recall that we are in a coordinate system in which (2.4) holds, with 
$Y^\#_i(p_0)=0,\ p_0=0$.  For $z$ close to $0$ in $\RR^{n-k}$ and $|t|<
\delta$, we define
$$
G(t,z)=y(t,0,z)=\Cal{F}^t (0,z),
\tag{2.17}
$$
where we set
$$
\Cal{F}^t=\Cal{F}^{t_1}_{Y_1}\circ\cdots\circ \Cal{F}^{t_k}_{Y_k}.
\tag{2.18}
$$

\proclaim{Proposition 2.3} There is a neighborhood $U_0$ of $(0,0)\in
\RR^k\times\RR^{n-k}$ and a neighborhood $U_1$ of $p_0\in\Cal{O}$ such that
$$
G:U_0\longrightarrow U_1
\tag{2.19}
$$
is a Lipschitz homeomorphism, with Lipschitz inverse.
\endproclaim
\demo{Proof} We want to show that if $(t,z)$ and $(s,w)$ are distinct points
in a small neighborhood of $(0,0)$, then $x_1=G(t,z)$ and $x_2=G(s,w)$
are not too close.  Note that
$$
\Cal{F}^{-t}(x_1)=(0,z),\quad \Cal{F}^{-t}(x_2)=\Cal{F}^{s-t}(0,w).
\tag{2.20}
$$
Since $\Cal{F}^{-t}$ is Lipschitz, we have
$$
|\Cal{F}^{-t}(x_1)-\Cal{F}^{-t}(x_2)|\le C|x_1-x_2|.
\tag{2.21}
$$
Meanwhile, since the span of
$Y_1,\dots,Y_k$ is transversal to $\{(0,z)\}$ near $(0,z)=(0,0)$, we have
$$
|(0,z)-\Cal{F}^{s-t}(0,w)|\ge C\Bigl(|z-w|+|s-t|\Bigr).
\tag{2.22}
$$
Comparing (2.20)--(2.22) yields
$$
|x_1-x_2|\ge C\Bigl(|z-w|+|s-t|\Bigr),
\tag{2.23}
$$
as desired.
\enddemo

$$\text{}$$
{\bf 3. The pull-back of a Levi-flat CR structure}
\newline {}\newline

In \S{2} we constructed a bi-Lipschitz map
$$
G:U_0\longrightarrow U_1,\quad G(t,z)=\Cal{F}^t(0,z),
\tag{3.1}
$$
taking sets $z=z_0$ to leaves of the foliation whose tangent space is the
involutive Lipschitz bundle $\Cal{E}\subset TU_1$.  Let us denote by
$\Cal{E}_0\subset TU_0$ the pull-back of $\Cal{E}$, so $\Cal{E}_0$ is
spanned by $\pa/\pa t_j,\ 1\le j\le k$.  Now we take $k=2m$ and suppose 
there is a complex structure on $\Cal{E},\ J\in \End(\Cal{E})$.  We pull 
this back to a complex structure $J^0\in \End(\Cal{E}_0)$, examine its
regularity, and show that if $J$ is formally integrable then so is $J^0$.

Since Lipschitz sections of $\Cal{E}$ are given as linear combinations 
over $\Lip(U_1)$ of the vector fields $Y_1,\dots,Y_k$, the action of $J$
is given by
$$
JY_i=\sum\limits_{j=1}^k \Cal{J}_{ij}(x)\, Y_j.
\tag{3.2}
$$
We can make various hypotheses on the regularity of $J$.  For example, we 
might assume
$$
\Cal{J}_{ij}\in \Lip(U_1),
\tag{3.3}
$$
or we might make the weaker hypothesis
$$
\Cal{J}_{ij}\in C^r(U_1),
\tag{3.4}
$$
for some $r\in (1/2,1)$.  In any case, the complex structure induced on 
$\Cal{E}_0$ is given by 
$$
J^0\, \frac{\pa}{\pa t_i}=\sum\limits_{j=1}^k \Cal{J}^0_{ij}(t,z)\, 
\frac{\pa}{\pa t_j},\quad \Cal{J}^0_{ij}(t,z)=\Cal{J}_{ij}(G(t,z)).
\tag{3.5}
$$
It is clear that
$$
\aligned
\Cal{J}_{ij}\in \Lip(U_1)&\Longrightarrow \Cal{J}^0_{ij}\in \Lip(U_0), \\
\Cal{J}_{ij}\in C^r(U_1)&\Longrightarrow \Cal{J}^0_{ij}\in C^r(U_0),
\endaligned
\tag{3.6}
$$
the latter provided $0<r<1$.  

We next discuss integrability conditions.  One approach would be to 
form the ``Nijenhuis'' tensor, associated to $J$ by
$$
N(X,Y)=[X,Y]-[JX,JY]+J[X,JY]+J[JX,Y],
\tag{3.7}
$$
for Lipschitz sections $X$ and $Y$ of $\Cal{E}$.  If $J$ is Lipschitz, then
(3.7) belongs to $L^\infty(U_1)$.  If $J$ satisfies (3.4) with $r>1/2$, 
then by Lemma 1.2 of [HT], the right side of (3.7) is a distribution,
belonging to $C^{r-1}_*(U_1)$.  Now such a singular distribution does not
necessarily pull back well under a bi-Lipschitz map.  Instead, we will
work on individual leaves.  

We start by defining $\Cal{N}^0_{z_0}$, associated with $J^0$, on a leaf 
in $(t,z)$-space where $z=z_0$ is constant.  We set
$$
\Cal{N}^0_{z_0}(X,Y)=[X,Y]-[J^0 X,J^0 Y]+J^0[X,J^0 Y]+J^0[J^0 X,Y],
\tag{3.8}
$$
where $X$ and $Y$ are linear combinations of $\pa/\pa t_i,\ 1\le i\le k$, 
and $J^0=J^0(t,z_0)$.  For each fixed $z_0$, this defines an element of 
$L^\infty(\Cal{O}_0)$ if $J^0$ is Lipschitz and an element of $C^{r-1}_*
(\Cal{O}_0)$ if (3.4) holds, with $r\in (1/2,1)$.  Here $\Cal{O}_0=\{
t\in\RR^k:|t|<\delta\}$.  As we have seen, for each $z_0$ close to $0$,
$G_{z_0}(t)=G(t,z_0)$ yields a $C^{1,1}$-diffeomorphism of $\Cal{O}_0$
onto a neighborhood of $x_0=G(0,z_0)$ in the leaf through $x_0$.  In light
of this, the following is useful.

\proclaim{Proposition 3.1} Assume $\varphi:\Cal{O}_0\rightarrow \Cal{O}_1$
is a $C^{1,1}$-diffeomorphism between open sets in $\RR^k$.  Then the 
pull-back
$$
\varphi^*:\Lip(\Cal{O}_1)\rightarrow \Lip(\Cal{O}_0),\quad
\varphi^*f(x)=f(\varphi(x))
\tag{3.9}
$$
extends to
$$
\varphi^*:H^{s,p}(\Cal{O}_1)\longrightarrow H^{s,p}(\Cal{O}_0),
\tag{3.10}
$$
for each $s\in [-1,1],\ p\in (1,\infty)$.  Furthermore, for each $r\in(0,1)$,
$$
\varphi^*:C^{r-1}_*(\Cal{O}_1)\longrightarrow C^{r-1}_*(\Cal{O}_0).
\tag{3.11}
$$
\endproclaim
\demo{Proof} The result (3.10) is easy for $s=0,1$, and follows by 
interpolation for $s\in (0,1)$.  Now suppose $s\in [-1,0)$.  We have, 
for compactly supported $u$,
$$
\aligned
(u,\varphi^*v)&=\int u(x) v(\varphi(x))\, dx \\
&=\int u(\varphi^{-1}(x))v(x)\, |\det D\varphi^{-1}(x)|\, dx.
\endaligned
\tag{3.12}
$$
We have $|\det D\varphi^{-1}|\in \Lip(\Cal{O}_1)$, hence (by the case already
treated) $u\in H^{\sigma,q}\Rightarrow (u\circ\varphi)|\det D\varphi^{-1}|
\in H^{\sigma,q}$, for $\sigma\in [0,1], q\in (1,\infty)$.  Thus by 
duality we have $v\in H^{s,p}\Rightarrow \varphi^*v\in H^{s,p}$, for
$s\in [-1,0), p\in (1,\infty)$, as desired.  Next, note that if $g\in \Lip
(\Cal{O}_1)$ and $X$ is a Lipschitz vector field on $\Cal{O}_1$, then
$\varphi$ transforms $X$ to a Lipschitz vector field $\widetilde{X}$ on 
$\Cal{O}_0$ and
$$
\varphi^*(Xg)=\widetilde{X}\, \varphi^*g,
\tag{3.13}
$$
as elements of $L^\infty(\Cal{O}_0)$.  Now if $g\in C^r(\Cal{O}_1)$, we have 
$\varphi^* g\in C^r(\Cal{O}_0)$ and then $\widetilde{X}\, \varphi^*g\in
C^{r-1}_*(\Cal{O}_0)$, which yields (3.11).
\enddemo

$\text{}$ \newline
{\smc Remark}.  More generally, if $\varphi$ is a diffeomorphism of class
$C^{1+r},\ r\in (0,1)$, then (3.13) holds with $\widetilde{X}$ a 
$C^r$-vector field.  Also, by Lemma 1.2 of [HT], 
$$
\tilde{g}\in C^r\Longrightarrow \widetilde{X} \tilde{g}\in C^{r-1}_*,\quad
\text{provided }\ r>\frac{1}{2},
\tag{3.14}
$$
so (3.11) still holds, as long as $r>1/2$.

$$\text{}$$
{\bf 4. The Newlander-Nirenberg theorem with parameters}
\newline {}\newline

The Newlander-Nirenberg theorem provides local holomorphic coordinates
on a manifold $\Omega$ with an almost complex structure satisfying the
formal integrability condition that its Nijenhuis tensor vanishes.  In the
setting of a relatively smooth almost complex structure $J$ the smooth
dependence of such coordinate functions on $J$ was noted in [NN] and played
a role in [Ni].  Here we aim to examine the dependence of such
coordinates on $J$, in appropriate function spaces, in the context of the
lower regularity hypotheses made here.  Verifying this regularity will
involve giving a review of the method of construction of holomorphic 
coordinates introduced in [M], with modifications as in [HT] to handle the
still weaker regularity hypotheses made here.

Given $p_0\in\Omega$, take coordinates $x=(x_1,\dots,x_{2m})$, centered at
$p_0$, with respect to which
$$
J(p_0) \frac{\pa}{\pa x_j}=\frac{\pa}{\pa x_{j+m}},\quad
J(p_0)\frac{\pa}{\pa x_{j+m}}=-\frac{\pa}{\pa x_j},\quad 1\le j\le m.
\tag{4.1}
$$
The condition for a function $f$, defined near $p_0$, to be holomorphic, 
is that $f$ be annihilated by the vector fields
$$
X_j=\frac{1}{2}\Bigl(\frac{\pa}{\pa x_j}+iJ\frac{\pa}{\pa x_j}\Bigr),\quad
1\le j\le m,
\tag{4.2}
$$
and in light of (4.1) we have $J(\pa/\pa x_j)=\pa/\pa x_{j+m}+
\sum_{\ell=1}^{2m}c_{j\ell} \pa/\pa x_\ell$ with $c_{j\ell}(0)=0\ (p_0=0)$.
Setting $y_j=x_{j+m},\ \pa/\pa z_j=(1/2)(\pa/\pa x_j-i\pa/\pa y_j),\
\pa/\pa \zbar_j=(1/2)(\pa/\pa x_j+i\pa/\pa y_j)$, we can write these 
complex vector fields as
$$
\frac{\pa}{\pa \zbar_j}+\sum\limits_{\ell=1}^m \Bigl(\alpha_{j\ell}
\frac{\pa}{\pa\zbar_\ell}+\beta_{j\ell}\frac{\pa}{\pa z_\ell}\Bigr),
\quad 1\le j\le m.
\tag{4.3}
$$
Next, by a device similar to that used in (2.2)--(2.4), we can take linear
combinations of these vector fields to obtain
$$
Z_j=\frac{\pa}{\pa\zbar_j}-\sum\limits_{\ell=1}^m a_{j\ell} 
\frac{\pa}{\pa z_\ell},\quad 1\le j\le m.
\tag{4.4}
$$
If $J$ is of class $C^r$, then the coefficients in (4.3) and (4.4) are also
of class $C^r$.

The formal integrability condition is that the Lie brackets $[X_j,X_\ell]$
are all linear combinations of $X_1,\dots,X_m$.  If $J\in C^1$, then 
$[X_j,X_\ell]$ is a linear combination with continuous coefficients.  If 
$J\in C^r$ with $r>1/2$, then the Lie brackets are still well defined, 
and the coefficients are distributions of class $C^{r-1}_*$.  In such a case,
it follows that the brackets $[Z_j,Z_\ell]$ are linear combinations of
$Z_1,\dots,Z_m$, which forces
$$
[Z_j,Z_\ell]=0,\quad 1\le j,\ell\le m.
\tag{4.5}
$$

It is convenient to use matrix notation.  Set $A_j=(a_{j1},\dots,a_{jm})$
(a row vector), $A=(a_{j\ell}),\ F=(f_1,\dots,f_m)$ (a row vector), and
$\pa/\pa\zbar=(\pa/\pa \zbar_1,\dots,\pa/\pa\zbar_m)^t$ (a column vector).
The condition that $f_1,\dots,f_m$ be $J$-holomorphic is that
$$
\frac{\pa F}{\pa \zbar}=A\, \frac{\pa F}{\pa z},
\tag{4.6}
$$
and the formal integrability condition (4.5) is
$$
\frac{\pa A_j}{\pa\zbar_\ell}+A_j \frac{\pa A_\ell}{\pa z}=
\frac{\pa A_\ell}{\pa \zbar_j}+A_\ell \frac{\pa A_j}{\pa z},\quad
1\le j,\ell\le m.
\tag{4.7}
$$
The proof of the Newlander-Nirenberg theorem consists of the construction
of $F$, mapping a neighborhood of $p_0$ in $\Omega$ diffeomorphically onto
a neighborhood of $0$ in $\CC^m$, and solving (4.6).

Malgrange's method constructs $F$ as a composition
$$
F=G\circ H.
\tag{4.8}
$$
Different techniques are applied to construct the diffeomorphisms $G$ and $H$.
We run through these constructions, paying particular attention to the 
dependence on the matrix $A$.  The Cauchy-Riemann equations (4.6) transform to
$$
\frac{\pa G}{\pa\zetabar}=B\, \frac{\pa G}{\pa\zeta},
\tag{4.9}
$$
for $\zeta=H(z)$, where $B$ is given by
$$
\frac{\pa H}{\pa\zbar}+\frac{\pa\overline{H}}{\pa\zbar}(B\circ H)=
A\Bigl[ \frac{\pa H}{\pa z}+\frac{\pa\overline{H}}{\pa z}(B\circ H)\Bigr],
\tag{4.10}
$$
or equivalently
$$
B\circ H=-\Bigl(\frac{\pa\overline{H}}{\pa \zbar}-
A\frac{\pa\overline{H}}{\pa z}\Bigr)^{-1}\Bigl(\frac{\pa H}{\pa\zbar}
-A\frac{\pa H}{\pa z}\Bigr).
\tag{4.11}
$$
The formal integrability condition (4.7) implies the corresponding
formal integrability of the new Cauchy-Riemann equations, i.e.,
$$
\frac{\pa B_j}{\pa\zetabar_\ell}+B_j\frac{\pa B_\ell}{\pa\zeta}
=\frac{\pa B_\ell}{\pa\zetabar_j}+B_\ell \frac{\pa B_j}{\pa\zeta},
\quad 1\le j,\ell\le m,
\tag{4.12}
$$
where $B_j$ are the rows of $B,\ 1\le j\le m$.  Furthermore, if $B$
satisfies (4.11), then the actual integrability, i.e., the existence of a 
diffeomorphism $G$ satisfying (4.9), is equivalent to the actual
integrability of $J$, i.e., the existence of a diffeomorphism $F$ 
satisfying (4.6).

A key idea of [M] to guarantee the existence of a diffeomorphism $G$ 
satisfying (4.9) is to construct $H$ in such a fashion that if $B$ is
defined by (4.11) then
$$
\sum\limits_j \frac{\pa B_j}{\pa\zeta_j}=0.
\tag{4.13}
$$
Equivalently, the task is to construct a diffeomorphism $H$ on a neighborhood
$U$ of $p_0=0$ in $\CC^m$ such that, if $B$ is defined by (4.11), then (4.13)
holds.  It is convenient to dilate the $z$-variable, so that $A(z)$ in (4.11)
is replaced by $A_t(z)=A(tz)$, and we solve on the unit ball, which we denote 
$U$, for sufficiently small positive $t$.  Note that if $A\in C^r_*$ and $A(0)
=0$, then $\|A_t\|_{C^r_*(\Ubar)}\rightarrow 0$ as $t\rightarrow 0$.  If we
relabel $A_t$ as $A$, we want to establish the following variant of Lemma 3.2
of [HT].  To state it, let us set
$$
\Cal{A}^r(\eta)=\Bigl\{A\in C^r_*(\Ubar):A(0)=0,\ \|A\|_{C^r(\Ubar)}<\eta
\Bigr\}.
\tag{4.14}
$$

\proclaim{Proposition 4.1}  Assume $r>1/2$.  Given $\ep,\delta>0$, 
there exists $\eta>0$ such that for any $A\in\Cal{A}^r(\eta)$ one can find
$$
H\in C^{1+r}_*(\Ubar),
\tag{4.15}
$$
satisfying
$$
H(0)=0,\quad \|H-\text{id}\|_{C^{1+r}_*(\Ubar)}<\delta,
\tag{4.16}
$$
and such that $B\in C^r_*(\Ubar)$, defined by (4.11), satisfies (4.13),
and $\|B\|_{L^\infty(U)}<\ep$.  Furthermore, $H$ is obtained as a $C^1$ map
$$
\Cal{A}^r(\eta)\longrightarrow C^{1+r}_*(\Ubar),\quad A\mapsto H,
\tag{4.17}
$$
\endproclaim
\demo{Proof}  Let us set
$$
\Phi(H,A)=E=-\Bigl(\frac{\pa \overline{H}}{\pa\zbar}-A 
\frac{\pa\overline{H}}{\pa z}\Bigr)^{-1} 
\Bigl(\frac{\pa H}{\pa\zbar}-A\frac{\pa H}{\pa z}\Bigr).
\tag{4.18}
$$
Then $\Phi$ is a $C^1$ map
$$
\Phi:\Cal{B}^{r+1}(\delta)\times\Cal{A}^{r}(1)\longrightarrow C^r(\Ubar),
\tag{4.19}
$$
where $\Cal{A}^r(\eta)$ is as in (4.14) and
$$
\Cal{B}^{r+1}(\delta)=\Bigl\{H\in C^{1+r}_*(\Ubar):H(0)=0,\
\|H-\text{id}\|_{C^{1+r}_*(\Ubar)}<\delta\Bigr\}.
\tag{4.20}
$$
If $B\circ H=\Phi(H,A)$, an application of the chain rule gives
$$
\frac{\pa B_j}{\pa \zeta_j}\circ H=\Bigl(\frac{\pa K}{\pa\zeta_j}\circ H
\Bigr)\frac{\pa E_j}{\pa z}+\Bigl(\frac{\pa\overline{K}}{\pa\zeta_j}\circ H
\Bigr)\frac{\pa E_j}{\pa\zbar},\quad K=H^{-1}.
\tag{4.21}
$$
Using the identity
$$
(DK)\circ H(z)=DH(z)^{-1}
\tag{4.22}
$$
of real $(2m)\times(2m)$ matrices, one can express $(\pa K/\pa\zeta_j)\circ H$
and $(\pa\overline{K}/\pa\zeta_j)\circ H$ in terms of the $z$- and 
$\zbar$-derivatives of $H$ and $\overline{H}$.  It follows from Lemma 1.2 of 
[HT] (extended to function spaces on bounded domains) that
$$
\Psi(H,A)=\sum\limits_j \frac{\pa B_j}{\pa\zeta_j}\circ H
\tag{4.23}
$$
defines a $C^1$ map
$$
\Psi:\Cal{B}^{r+1}(\delta)\times\Cal{A}^{r}(1)\longrightarrow 
C^{r-1}_*(\Ubar).
\tag{4.24}
$$
In fact $H\mapsto\Psi(H,A)$ is given by a nonlinear second order differential 
operator:
$$
\Psi(H,A)=\sum\limits_j a_j(\nabla H)\, \pa_j b_j(A,\nabla H),
\tag{4.25}
$$
where $a_j$ and $b_j$ are smooth in their arguments.  We note that if
$$
H(z)=z+\ep h(z),
\tag{4.26}
$$
then
$$
\Phi(H,0)=-\ep \frac{\pa h}{\pa\zbar}+O(\ep^2),
\tag{4.27}
$$
and (for $A=0$)
$$
\frac{\pa B_j}{\pa \zeta_j}\circ H=-\ep \frac{\pa^2 h}{\pa z_j\pa\zbar_j}
+O(\ep^2).
\tag{4.28}
$$
Hence
$$
\Psi(\text{id},0)=0,
\tag{4.29}
$$
and
$$
D_H \Psi(\text{id},0)h=-\sum\limits_j \frac{\pa^2 h}{\pa z_j\pa\zbar_j}
=-\frac{1}{4}\, \Delta h.
\tag{4.30}
$$
The map (4.30) has a right inverse
$$
\widetilde{G}h=-4(Gh-Gh(0)),
\tag{4.31}
$$
where $G$ denotes the solution operator to
$$
\Delta v=h\ \text{ on }\ U,\quad v\bigr|_{\pa U}=0,
\tag{4.32}
$$
which has the mapping property
$$
G:C^{r-1}_*(\Ubar)\longrightarrow C^{r+1}_*(\Ubar),
\tag{4.33}
$$
valid for $r>0$.  From here, Proposition 4.1 follows from the Implicit 
Function Theorem.
\enddemo

If $A\in\Cal{A}^r(\eta)$ satisfies the formal integrability condition (4.7) 
and we construct $H$ according to Proposition 4.1, defining $B$ by (4.11), 
then $B$ satisfies both (4.12) and (4.13).  This is an overdetermined
elliptic system (if $\ep$ is small enough), which we will write as
$$
\sum\limits_{|\alpha|=1} a_\alpha(B)\, \pa^\alpha B=0.
\tag{4.34}
$$
The a priori regularity we have on $B$ from (4.11) is
$$
B\circ H\in C^r(\Ubar),\quad \text{hence }\ B\in C^r(\Obar),
\tag{4.35}
$$
where $\Obar\subset H^{-1}(\Ubar)$.  As shown in Lemma 4.1 of [HT], having
this a priori information with $r>1/2$ allows us to obtain
$$
B\in C^N_{\text{loc}}(\Cal{O}),
\tag{4.36}
$$
for each $N<\infty$.  Then classical results yield
$$
|\pa^\alpha B(\zeta)|\le C^{|\alpha|+1} \alpha!,\quad \zeta\in\Cal{O}^b
\subset\subset \Cal{O},\quad C=C(\Cal{O}^b).
\tag{4.37}
$$
Once we have this (as [M] noted), producing a diffeomorphism $G$ such that 
(4.9) holds, which amounts to proving the Newlander-Nirenberg theorem in
the real analytic setting, is amenable to classical techniques for solving
real analytic systems of partial differential equations.  A self contained
treatment of a complex Frobenius theorem in the real analytic category,
which will produce such a construction, is presented in Appendix A of this
paper.

Having described how to obtain the holomorphic coordinate system (4.8),
we want to examine how it depends on $A$.  So we pick
$$
A_1, A_2\in\Cal{A}^r(\eta),
\tag{4.38}
$$
with $r>1/2$ and $\eta>0$ sufficiently small, and turn to the task of 
estimating, in turn (with obvious notation), $H_1-H_2,\ B_1-B_2,\
G_1-G_2$, and then $F_1-F_2$, in terms of $A_1-A_2$.  The assertion from 
Proposition 4.1 that the map (4.17) is $C^1$ leads immediately to our 
first estimate:
$$
\|H_1-H_2\|_{C^{1+r}_*(\Ubar)}\le C \|A_1-A_2\|_{C^r_*(\Ubar)}.
\tag{4.39}
$$
We also have a $C^1$ map
$$
\Cal{A}^r(\eta)\rightarrow C^r_*(\Ubar),\quad A\mapsto \wtB=B\circ H,
\tag{4.40}
$$
in light of the formula (4.11).  Hence
$$
\|\wtB_1-\wtB_2\|_{C^r_*(\Ubar)}\le C \|A_1-A_2\|_{C^r_*(\Ubar)}.
\tag{4.41}
$$

Now $B_j=\wtB_j\circ K_j$ with $K_j=H_j^{-1}$.  While $A\mapsto H$ is $C^1$
from $\Cal{A}^r(\eta)$ to $C^{r+1}_*(\Ubar)$, one has that $A\mapsto K=
H^{-1}$ is a continuous map from $\Cal{A}^r(\eta)$ to $C^{r+1}_*(\Obar)$
and a $C^1$ map to $C^r_*(\Obar)$, where $\Obar$ is a neighborhood of $0$
containing $H^{-1}(\Ubar)$ for all $H$ as in (4.16).  Consequently
$$
\|K_1-K_2\|_{C^r_*(\Obar)}\le C\|A_1-A_2\|_{C^r_*(\Ubar)},\quad
\|K_j\|_{C^{1+r}_*(\Obar)}\le C.
\tag{4.42}
$$
Let us write
$$
B_1-B_2=\wtB_1\circ K_1-\wtB_2\circ K_1+\wtB_2\circ K_1-\wtB_2\circ K_2.
\tag{4.43}
$$
We have
$$
\|\wtB_1\circ K_1-\wtB_2\circ K_1\|_{C^r(\Obar)}\le C\|\wtB_1-\wtB_2
\|_{C^r(\Ubar)} \|K_1\|_{C^1(\Obar)},\quad 0<r\le 1,
\tag{4.44}
$$
and
$$
\|\wtB_2\circ K_1-\wtB_2\circ K_2\|_{L^\infty(\Obar)}\le 
\|\wtB_2\|_{C^r(\Ubar)} \|K_1-K_2\|^r_{L^\infty(\Obar)},\quad 0<r\le 1.
\tag{4.45}
$$
Putting together (4.43)--(4.45), using the estimates (4.41)--(4.42),
we obtain
$$
\|B_1-B_2\|_{L^\infty(\Obar)}\le C(\|A_2\|_{C^\rho(\Ubar)})
\|A_1-A_2\|^{\rho}_{C^s(\Ubar)},\quad
\frac{1}{2}<\rho,s<1.
\tag{4.46}
$$
(It is convenient to replace $r$ by $\rho$ in our use of (4.45) and to
replace $r$ by $s$ in our use of (4.42) and (4.44).  Typically we will want 
to take $\rho$ as large as possible and $s$ as small as possible.)
The estimate (4.46) is a relatively weak estimate, 
a consequence of the rather rough
dependence of $\wtB\circ K$ on $K$.  Fortunately, (4.46) can be improved
substantially via use of the fact that $B_1$ and $B_2$ both satisfy the 
elliptic system (4.34).  Hence $V=B_1-B_2$ solves
$$
\sum\limits_{|\alpha|=1} a_\alpha(B_1)\, \pa^\alpha V=
\sum\limits_{|\alpha|=1} [a_\alpha(B_2)-a_\alpha(B_1)]\, \pa^\alpha B_2.
\tag{4.47}
$$
In fact, as one sees from (4.12)--(4.13), $a_\alpha(B)=a_\alpha^0+M_\alpha 
B$, with $M_\alpha$ a linear map, and hence $V$ solves the linear elliptic 
system (with real analytic coefficients)
$$
\sum\limits_{|\alpha|=1} a_\alpha(B_1)\, \pa^\alpha V
-\sum\limits_{|\alpha|=1} (\pa^\alpha B_2)\, M_\alpha V=0.
\tag{4.48}
$$
The estimates (4.37) hold for $B_1$ and $B_2$.  Local elliptic regularity
results yield
$$
\bigl|\pa^\alpha\bigl(B_1(\zeta)-B_2(\zeta)\bigr)\bigr|\le
C^{|\alpha|+1} \alpha!\, \|B_1-B_2\|_{L^\infty(\Obar)},\quad
\zeta\in\Cal{O}^b\subset\subset \Cal{O}.
\tag{4.49}
$$
Then the method of solving (4.9) covered in Appendix A gives
$$
\bigl|\pa^\alpha \bigl(G_1(\zeta)-G_2(\zeta)\bigr)\bigr|\le
C^{|\alpha|+1} \alpha!\, \|B_1-B_2\|_{L^\infty(\Obar)},\quad
\zeta\in\Cal{O}^b.
\tag{4.50}
$$

Now, with $F_j=G_j\circ H_j$, we can set
$$
F_1-F_2=G_1\circ H_1-G_2\circ H_1+G_2\circ H_1-G_2\circ H_2.
\tag{4.51}
$$

Under the bounds on $H_j$ in $C^{1+r}$ and on $G_j$ in $C^N$ produced 
above, we have, for $U^b\subset\subset U,\ r\in (1/2,1)$,
$$
\aligned
\|G_1\circ H_1-G_2\circ H_1\|_{C^{1+r}(U^b)}&\le C\|G_1-G_2
\|_{C^{1+r}(\Cal{O}^b)} \\ 
&\le C\|B_1-B_2\|_{L^\infty(\Obar)},
\endaligned
\tag{4.52}
$$
and
$$
\aligned
\|G_2\circ H_1-G_2\circ H_2\|_{C^{1+r}(U^b)}&\le C\|G_2\|_{C^3(\Cal{O}^b)}
\|H_1-H_2\|_{C^{1+r}(\Ubar)} \\ &\le C\|A_1-A_2\|_{C^r(\Ubar)}.
\endaligned
\tag{4.53}
$$
Hence
$$
\aligned
\|F_1-F_2\|_{C^{1+r}(U^b)}&\le C\Bigl(\|B_1-B_2\|_{L^\infty(\Obar)}
+\|A_1-A_2\|_{C^r(\Ubar)} \Bigr) \\
&\le C(\|A_2\|_{C^\rho(\Ubar)}) \|A_1-A_2\|^\rho_{C^s(\Ubar)}+
C\|A_1-A_2\|_{C^r(\Ubar)},
\endaligned
\tag{4.54}
$$
given $1/2<r,s,\rho<1$.  For the last inequality, we have used (4.46).
As in that estimate, we typically want to
take $\rho$ as large as possible and $s$ as small as possible.

$$\text{}$$
{\bf 5. Structure of Levi-flat CR-manifolds}
\newline {}\newline

In this section we assume $\Cal{S}$ is a Lipschitz subbundle of $\CC T\Omega$,
satisfying
$$
\Cal{S}_p\cap \Sbar_p=0,\quad \forall\ p\in\Omega.
\tag{5.1}
$$
Hence $\Cal{S}_p+\Sbar_p$ has constant dimension (say $k$), and so does
$\Cal{E}_p$, defined by (1.6).  It follows that $\Cal{E}$ and $\Cal{S}+\Sbar$ 
are Lipschitz vector bundles, and of course $\Cal{V}=0$.  The bundle $\Cal{E}
\subset T\Omega$ gets a complex structure
$$
J\in\Lip(\Omega,\End \Cal{E}),
\tag{5.2}
$$
and
$$
\Cal{S}_p=\{u+iJu:u\in\Cal{E}_p\}.
\tag{5.3}
$$
We make the involutivity hypotheses (1.1)--(1.2).  As explained in the 
introduction, this is equivalent to the hypothesis that $\Cal{E}$ is
involutive plus the hypothesis that the Nijenhuis tensor of $J$ vanishes. 
A manifold $\Omega$ with such a structure $(\Cal{E},J)$ is said to be a 
Levi-flat CR-manifold.

In this setting, a function $f$ on an open set $\Cal{O}\subset\Omega$ is
called a CR function provided
$$
Zf=0\ \text{ on }\ \Cal{O},\quad \forall\ Z\in\Lip(\Cal{O},\Cal{S}),
\tag{5.4}
$$
or equivalently
$$
Xf+i(JX)f=0\ \text{ on }\ \Cal{O},\quad \forall\ X\in\Lip(\Cal{O},\Cal{E}).
\tag{5.5}
$$
Given the regularity of $X$ and $Z$, we see that $Zf$ is a well defined
distribution for any $f\in L^2_{\text{loc}}(\Cal{O})$.  Our goal here is to
construct a rich class of CR functions $f$ having the regularity
$$
f,\ Xf\in C^s(\Cal{O}),\quad \forall\ X\in \Lip(\Cal{O},\Cal{E}),
\tag{5.6}
$$
given $s<1/2$.  In fact $f$ and $Xf$ will have further regularity along the 
leaves of the foliation tangent to $\Cal{E}$, as will be explained below.

To begin the construction of such CR functions, we implement the results
of \S\S{2--3}.  For any $p\in\Omega$, there are a neighborhood $U_1$ of $p$,
a neighborhood $U_0$ of $0\in \RR^n$ ($n=\text{dim}\, \Omega$) and a 
bi-Lipschitz map $\Cal{G}:U_0\rightarrow U_1$, 
pulling $\Cal{E}$ back to the bundle
$\Cal{E}^\#$ spanned by $\pa/\pa t_1,\dots \pa/\pa t_k$, where in $U_0\subset
\RR^k\times\RR^{n-k}$ we have coordinates $(t,z)=(t_1,\dots,t_k,z_1,\dots,
z_{n-k})$.  Furthermore, Lipschitz sections of $\Cal{E}$ are transformed to
Lipschitz vector fields on $U_0$, and $J$ is transformed to
$$
J_0\in\Lip(U_0,\End \Cal{E}^\#).
\tag{5.7}
$$
We may as well assume $U_0=U'_0\times U''_0$, where $U'_0$ is a neighborhood
of $0\in\RR^k$ and $U''_0$ a neighborhood of $0\in \RR^{n-k}$.  Then $J_0=
J_0(z)$ is effectively a family of integrable almost complex structures on 
$U'_0$, parametrized by $z\in U''_0$.  Of course $k$ is even; say $k=2m$.

Now we can apply the results of \S{4}.  We construct holomorphic functions 
$F=(f_1,\dots,f_m)$ on $U'_0$, depending on $z$ as a parameter, say $F=F_z:
U'_0\rightarrow \CC^m,\ z\in U''_0$.  (Note that $z$ has a different role 
here than in \S{4}; this should not cause confusion.)  We construct $F_z$
as a composition:
$$
F_z(t)=G_z(H_z(t)).
\tag{5.8}
$$
The family of diffeomorphisms $H_z$ is constructed in Proposition 4.1, via
an implicit function theorem.  Perhaps after shrinking $U'_0$ and $U''_0$, 
we have $H_z\in C^{1+r}(U'_0)$ for each $z\in U'_0$, given $r<1$, and
$$
\aligned
\|H_z-H_{z'}\|_{C^{1+r}(U'_0)}&\le C\|A_z-A_{z'}\|_{C^2(U'_0)} \\
&\le C |z-z'|^{1-r},
\endaligned
\tag{5.9}
$$
if $1/2<r<1$.  Here we have used
$$
\|A_1-A_2\|_{C^r}\le C\|A_1-A_2\|_{L^\infty}^{1-r} 
\|A_1-A_2\|_{\Lip}^r,\quad 0<r<1.
\tag{5.10}
$$
As explained in \S{4}, the construction of $G_z$ follows from the 
real-analytic version of the Newlander-Nirenberg theorem, a presentation of 
which is given here, in Appendix A.  Then we obtain $F_z=G_z\circ H_z$, and, 
by (4.54), with $U_0^b\subset\subset U'_0$, 
$$
\|F_z-F_{z'}\|_{C^{1+r}(U^b_0)}\le C(\|A_{z'}\|_{C^\rho(U'_0)})
\|A_z-A_{z'}\|_{C^s(U'_0)}^\rho+C\|A_z-A_{z'}\|_{C^r(U'_0)},
\tag{5.11}
$$
given $1/2<r,s,\rho<1$.  Here we pick $\rho=1-\ep,\ s=1/2+\ep$, and use 
(5.10) to obtain
$$
\aligned
\|F_z-F_{z'}\|_{C^{1+r}(U^b_0)}&\le C |z-z'|^{1/2-\delta}+C|z-z'|^{1-r} \\
&\le C|z-z'|^{1-r},
\endaligned
\tag{5.12}
$$
given $r\in (1/2,1)$, and taking $\ep$ (hence $\delta$) sufficiently small.

The functions $f_j(t,z)$ given by $F_z(t)=(f_1(t,z),\dots,f_m(t,z))$ are 
CR functions on $U_0$.  In addition, the functions $\varphi_j(t,z)=z_j,\
1\le j\le n-k$, are CR functions on $U_0$.  Then
$$
\Phi(t,z)=(f_1(t,z),\dots,f_m(t,z),z_1,\dots,z_{n-k})
\tag{5.13}
$$
gives a H{\"o}lder continuous homeomorphism of $U_0$ (possibly shrunken 
some more) onto an open subset of $\CC^m\times\RR^{n-k}$.  We compose with 
$G^{-1}$ to get associated CR functions on $U_1\subset\Omega$.  Let us 
formally record the result.

\proclaim{Proposition 5.1} Given $\Omega$ with a Lipschitz, Levi-flat CR 
structure, $p\in\Omega$, there exists a neighborhood $U_1$ of $p$ and a 
homeomorphism
$$
\Phi:U_1\longrightarrow \Cal{O}\subset\CC^m\times\RR^{n-k},
\tag{5.14}
$$
whose components are CR functions $\varphi_1,\dots,\varphi_{m+n-k}$ on
$U_1$.  We have
$$
\varphi_j,\ X\varphi_j\in C^s(U_1),\quad \forall\ X\in\Lip(U_1,\Cal{E}),
\tag{5.15}
$$
for any $s<1/2$.  Furthermore, $\Phi$ is a $C^{1+r}$-embedding of each leaf
in $U_1$, tangent to $\Cal{E}$, into $\CC^m\times\RR^{n-k}$, for each $r<1$.
\endproclaim

$\text{}$ \newline
{\smc Remark}.  Note that if $\psi$ is a smooth function on a neighborhood of 
the range of $\Phi$ in $\CC^m\times\RR^{n-k}$ and if $\psi$ is holomorphic
in the $\CC^m$-variables, then $\psi(\varphi_1,\dots,\varphi_{m+n-k})$
is a CR function on $U_1$.
\newline $\text{}$

If $\text{dim}\, \Cal{S}_p=1$, so $k=2$ and the leaves tangent to $\Cal{E}$
are 2-dimensional, then we can use the results of Appendix B in place of 
those of \S{4}.  Consequently we can improve the regularity result (5.15) to
$$
|\varphi_j(x)-\varphi_j(x')|,\ |X\varphi_j(x)-X\varphi_j(x')|\le C
\sigma^\#(|x-x'|),
\tag{5.16}
$$
where, given $a>0$,
$$
\sigma^\#(\delta)=\delta \Bigl(\log \frac{e}{\delta}\Bigr)^{1+a},
\tag{5.17}
$$
for $0<\delta\le 1$.

We now give a sufficient condition for the existence of a CR embedding $\Phi$
as in (5.14) that is a $C^1$ diffeomorphism.

\proclaim{Proposition 5.2} Assume $\Omega$ is a Levi-flat CR manifold with 
a CR structure regular of class $C^\rho$, with $\rho>3/2$. Then the map
$\Phi$ in (5.14) can be taken to be a $C^1$ diffeomorphism.
\endproclaim
\demo{Proof} The new regularity hypothesis is that $\Cal{S}$ is a $C^\rho$
bundle.  Thus $\Cal{E}$ and $J$ are regular of class $C^\rho$, and these
structures pull back to $C^\rho$ structures under the map $\Cal{G}$, which
is a $C^\rho$ diffeomorphism.  In particular, $A(t,z)$ is a $C^1$ function
of $z$ with values in $C^s(U'_0)$, with $s=\rho-1>1/2$.  Thus the implicit 
function theorem argument of Proposition 4.1 yields $H_z$, a $C^1$ function
of $z$ with values in $C^{1+s}$.  From here, one obtains $C^1$ dependence
of $G_z$ on $z$ and the result follows.
\enddemo

Note that if the leaves tangent to $\Cal{E}$ are 2-dimensional, we can 
obtain the conclusion of Proposition 5.2 whenever $\rho>1$, again using the 
results of Appendix B in place of those of \S{4}.

$\text{}$ \newline
{\smc Remarks on the embedded case}.
Suppose $\Omega\subset\CC^N$ is a $C^{1,1}$ submanifold, of real dimension 
$d$, and that $T_p\Omega\cap JT_p\Omega=\Cal{E}_p$ has constant real dimension 
$k=2m$, so $\Omega$ has the structure of a CR-manifold.  The vector bundle
$\Cal{E}\subset T\Omega$ is a Lipschitz vector bundle, and the condition
that $\Cal{E}$ be involutive is equivalent to the condition that $\Omega$
is a Levi-flat CR-manifold.  In such a case, the results of \S{2} imply
that $\Omega$ is foliated by manifolds, of real dimension $k$, tangent to
$\Cal{E}$, and smooth of class $C^{1,1}$.

In this case one does not need the Newlander-Nirenberg theorem (or a
refinement) to establish that these leaves are complex manifolds.  Rather
methods going back to Levi-Civita [LC], and developed further in [Som],
[Fr], and [Pin] suffice.  Levi-Civita's result for a single leaf is:

\proclaim{Proposition 5.3} If $M$ is a $C^1$ submanifold of $\CC^N$ and each 
tangent space $T_pM$ is $J$-invariant, then $M$ is a complex manifold.
\endproclaim
\demo{Proof} Fix $p\in M$, and represent $M$ near $p$ as the graph  over
the complex vector space $V=T_pM$; 
so one has a $C^1$ diffeomorphism $G:\Cal{O}
\rightarrow M$ where $\Cal{O}$ is a neighborhood of $0\in V$.  It is
readily verified that $DG(q)$ is $\CC$-linear for $q\in\Cal{O}$, so 
$G:\Cal{O}\rightarrow \CC^N$ is holomorphic.
\enddemo

In the setting above, we have a family $M_z$ of leaves, depending in a
Lipschitz fashion of $z\in\Cal{U}\subset\RR^\ell$, where $d=\ell+k$.
Given $p\in\Omega$, say $p\in M_{z_0}$, pick $V=T_p M_{z_0}$, and for $z$
close to $z_0$ we have $M_{z}$ locally a graph over $\Cal{O}\subset V$.
The comments above give local holomorphic diffeomorphisms $G_z:\Cal{O}
\rightarrow M_z\subset\CC^N$.  This construction, as we have said, is 
essentially classical.  The one point to make here is that we have the
Frobenius theory of [Ha], so we are able to treat submanifolds of class 
$C^{1,1}$ while previous treatments take $\Omega$ to be of class $C^2$.  
In connection with this, 
we note that Theorem 2.1 of [Pin] refers to CR-manifolds in
$\CC^N$ of class $C^m$, with $m\ge 1$, but a perusal of the proof shows that
the author means to say the relevant tangent spaces are smooth of class 
$C^m$, which holds if $\Omega\subset\CC^N$ is a submanifold of class $C^{m+1}$
(satisfying the CR property).

$$\text{}$$
{\bf 6. The complex Frobenius theorem}
\newline {}\newline

We recall our set-up.  We have a Lipschitz bundle $\Cal{S}\subset\CC 
T\Omega$, we assume $\Cal{S}+\Sbar$ is also a Lipschitz bundle, and we assume that
$$
X,Y\in\Lip(\Omega,\Cal{S})\Rightarrow [X,Y]\in L^\infty(\Omega,\Cal{S}),\ \
[X,\overline{Y}]\in L^\infty(\Cal{S}+\Sbar).
\tag{6.1}
$$
We then form the Lipschitz bundles $\Cal{V}\subset\Cal{E}\subset T\Omega$,
with fibers
$$
\Cal{V}_p=\Cal{S}_p\cap T_p\Omega,\quad \Cal{E}_p=\{w+\overline{w}:
w\in\Cal{S}_p\},
\tag{6.2}
$$
which therefore satisfy
$$
\aligned
X,Y\in\Lip(\Omega,\Cal{E})&\Longrightarrow [X,Y]\in L^\infty(\Omega,\Cal{E}),
\\
X,Y\in\Lip(\Omega,\Cal{V})&\Longrightarrow [X,Y]\in L^\infty(\Omega,\Cal{V}).
\endaligned
\tag{6.3}
$$
Furthermore, we have a complex structure on $\Cal{E}/\Cal{V}$,
$$
J\in \Lip(\Omega,\End \Cal{E}/\Cal{V}),
\tag{6.4}
$$
satisfying
$$
J(u\ \text{mod}\,\Cal{V})=v\ \text{mod}\, \Cal{V},\quad u+iv\in\Cal{S}_p.
\tag{6.5}
$$

Our proximate goal is to construct a Levi-flat CR manifold as a quotient
(locally) of $\Omega$, 
via the action of a local group of flows generated by sections
of $\Cal{V}$.  In order to achieve this, we need a further hypothesis on the 
regularity with which $\Cal{V}$ sits in $\Cal{E}$.  One way to put it is
the following.  Say $\text{dim}\,\Cal{V}_p=\ell\le k=\text{dim}\,\Cal{E}_p$.

\proclaim{Hypothesis V}  Each $p\in\Omega$ has a neighborhood $U_1$ 
on which there is a local Lipschitz frame field 
$\{X_1,\dots,X_k\}$ for $\Cal{E}$, such that
$\{X_1,\dots,X_\ell\}$ is a local frame field for $\Cal{V}$ and
$$
[X_i,X_j]=0,\quad 1\le i,j\le k.
\tag{6.6}
$$
\endproclaim

Later we will give other conditions that imply Hypothesis V, but for now
we show how it leads to the desired quotient space.

With respect to such a local frame field, for $x\in U_1$ we can identify 
$\Cal{E}_x/\Cal{V}_x$ with the linear span of $X_{\ell+1}(x),\dots,X_k(x)$,
and we can represent $J$ by a $(k-\ell)\times(k-\ell)$ matrix:
$$
JX_j=\sum\limits_{m=\ell+1}^k \Cal{J}_{jm}(x)\, X_m\ \text{mod}\, \Cal{V}_x,
\quad \ell+1\le j\le k.
\tag{6.7}
$$
Note that if $Y_j\in\Lip(U_1,\Cal{E})$ and $X_j+iY_j\in \Lip(U_1,\Cal{S})$,
so $Y_j=JX_j\ \text{mod}\, \Cal{V}$, we have
$$
\aligned
[X_i,X_j+iY_j]=i[X_i,Y_j]&\in L^\infty(\Omega,\Cal{S})\cap iL^\infty(\Omega,
\Cal{E}) \\ &\subset i L^\infty(\Omega,\Cal{V}),
\endaligned
\tag{6.8}
$$
for $1\le i\le\ell,\ \ell+1\le j\le k$, by (6.1) and (6.6).  Taking $Y_j$
to be the sum in (6.7), and noting that
$$
\Bigl[X_i,\sum\limits_{m=\ell+1}^k \Cal{J}_{jm}X_m\Bigr]
=\sum\limits_{m=\ell+1}^k (X_i\Cal{J}_{jm})\, X_m,
\tag{6.9}
$$
again by (6.6), we deduce that (6.9) actually vanishes, and hence
$$
X_i\Cal{J}_{jm}=0,\quad 1\le i\le\ell,\ \ell+1\le j,m\le k.
\tag{6.10}
$$

In a fashion parallel to (2.17) and (3.1), we set
$$
G(t,z)=\Cal{F}^t(0,z),\quad \Cal{F}^t=\Cal{F}^{t_1}_{X_1}\circ\cdots\circ
\Cal{F}^{t_k}_{X_k},
\tag{6.11}
$$
with $X_1,\dots,X_k$ as in Hypothesis V.  By Proposition 2.3, $G:U_0
\rightarrow U_1$ is a bi-Lipschitz map from a neighborhood $U_0$ of $(0,0)
\in \RR^k\times\RR^{n-k}$ to a neighborhood $U_1$ of $p\in\Omega$.
We denote by $\Cal{V}_0\subset TU_0$ the pull back of $\Cal{V}$, by
$\Cal{E}_0\subset TU_0$ the pull back of $\Cal{E}$, and by $\Cal{S}_0
\subset\CC TU_0$ the pull back of $\Cal{S}$.  Note that $\Cal{V}_0$ is
spanned by $\pa/\pa t_j,\ 1\le j\le\ell$ and $\Cal{E}_0$ by $\pa/\pa t_j,\
1\le j\le k$.  The quotient bundle $\Cal{E}_0/\Cal{V}_0$ is isomorphic to
the span of $\pa/\pa t_j$ for $\ell+1\le j\le k$, and the complex structure
$J$ on $\Cal{E}/\Cal{V}$ pulls back to $J^0$, given by
$$
J^0 \, \frac{\pa}{\pa t_j}=\sum\limits_{m=\ell+1}^k \Cal{J}^0_{jm}(t,z)\,
\frac{\pa}{\pa t_m},\quad \Cal{J}^0_{jm}(t,z)=\Cal{J}_{jm}(G(t,z)),\quad
\ell+1\le j\le k.
\tag{6.12}
$$
The result (6.10) is equivalent to
$$
\frac{\pa}{\pa t_i}\, \Cal{J}^0_{jm}(t,z)=0,\quad 1\le i\le \ell,\ 
\ell+1\le j,m\le k,
\tag{6.13}
$$
so we can write
$$
\Cal{J}^0_{jm}=\Cal{J}^0_{jm}(t'',z),\quad t''=(t_{\ell+1},\dots,t_k).
\tag{6.14}
$$

At this point it is natural to form the quotient space $\widetilde{U}_0=
U_0/\sim$, where we use the equivalence relation
$$
(t,z)\sim(s,z)\Longleftrightarrow (t_{\ell+1},\dots,t_k)=
(s_{\ell+1},\dots,s_k).
\tag{6.15}
$$
In other words, $\widetilde{U}_0$ is a neighborhood of $(0,0)\in \RR^{k-\ell}
\times\RR^{n-k}$, with coordinates
$$
(t_{\ell+1},\dots,t_k,z_1,\dots,z_{n-k}).
\tag{6.16}
$$
Note that $U_1$ fibers over $\widetilde{U}_0$, via
$$
\pi=P\circ G^{-1}:U_1\longrightarrow \widetilde{U}_0,
\tag{6.17}
$$
where $P(t_1,\dots,t_k)=(t_{\ell+1},\dots,t_k)$.
We will display a Levi-flat CR structure on $\widetilde{U}_0$, with
$\widetilde{\Cal{E}}_0$ the span of $\pa/\pa t_j,\ \ell+1\le j\le k$ and
$$
\widetilde{J}^0\, \frac{\pa}{\pa t_j}=\sum\limits_{m=\ell+1}^k
\Cal{J}^0_{jm}(t'',z)\, \frac{\pa}{\pa t_m}.
\tag{6.18}
$$
To see this, note that a vector field of the form
$$
\frac{\pa}{\pa t_j}+i\widetilde{J}^0\, \frac{\pa}{\pa t_j},\quad 
\ell+1\le j\le k,
\tag{6.19}
$$
can be regarded as a vector field on either $\widetilde{U}_0$ or $U_0$.
In the latter guise it is a Lipschitz section of $\Cal{S}_0$.  The
involutivity condition (6.1) has a counterpart for $\Cal{S}_0$, which
implies that the Nijenhuis tensor of $\widetilde{J}^0$ vanishes, so 
$\widetilde{U}_0$ has a Levi-flat CR structure, associated with 
$\widetilde{\Cal{S}}_0$, the span of vectors of the form (6.19).  This 
establishes the main result of this section, which we state formally.

\proclaim{Proposition 6.1} Assume $\Cal{S}$ and $\Cal{S}+\Sbar$ are
Lipschitz subbundles of $\CC T\Omega$, satisfying the involutivity 
condition (6.1) and also Hypothesis V.  Then each $p\in\Omega$ has a
neighborhood $U_1$ and a Lipschitz fibration $\pi:U_1\rightarrow 
\widetilde{U}_0$ onto a Levi-flat CR manifold, associated to a Lipschitz 
subbundle $\widetilde{\Cal{S}}_0\subset \CC T\widetilde{U}_0$, such that 
$$
\Cal{S}\Bigr|_{U_1}=(D\pi)^{-1}\, \widetilde{\Cal{S}}_0
\Bigr|_{\widetilde{U}_0}.
\tag{6.20}
$$
\endproclaim

We show that additional regularity conditions on $\Cal{V}$ and $\Cal{E}$
imply Hypothesis V.

\proclaim{Proposition 6.2} Assume each $p\in\Omega$ has a neighborhood on
which there is a frame field $\{W_1,\dots,W_k\}$ for $\Cal{E}$, of class
$C^{1,1}$, such that $\{W_1,\dots,W_\ell\}$ is a local frame field for 
$\Cal{V}$.  Then Hypothesis V holds.
\endproclaim
\demo{Proof} We begin with a construction parallel to (2.2)--(2.6), obtaining 
a local $C^{1,1}$ frame field $\{Y_1,\dots,Y_k\}$ for $\Cal{E}$ such that
$[Y_i,Y_j]=0$ for $1\le i,j\le k$ (though $\{Y_1,\dots,Y_\ell\}$ might not be
a local frame field for $\Cal{V}$).  As in (2.17)--(2.18), construct a 
diffeomorphism $G$, of class $C^{1,1}$, via which $Y_j$ are transformed to
$\pa/\pa t_j,\ 1\le j\le k$, and note that $W_j$ are transformed to $C^{1,1}$
vector fields $V_j=\sum_{i=1}^k v_{ji}(t,z)\, \pa/\pa t_i$.  Now produce
a $C^{1,1}$ diffeomorphism $H$ that straightens appropriate linear 
combinations of $V_1,\dots,V_\ell$ to $\pa/\pa s_1,\dots,\pa/\pa s_\ell$,
while each $\pa/\pa s_j\ (1\le j\le k)$ is a linear combination of 
$\pa/\pa t_1,\dots,\pa/\pa t_k$.  Then transform $\pa/\pa s_j$ via 
$(H\circ G)^{-1}$, to obtain the Lipschitz vector fields $X_j$ of Hypothesis 
V.
\enddemo

$$\text{}$$
{\bf A. A Frobenius theorem for real analytic, complex vector fields}
\newline {}\newline

Let $X_1,\dots,X_m$ be real analytic, complex vector fields on an open set
$\Cal{O}\subset\RR^n$.  We assume
$$
[X_k,X_\ell]=0,\quad 1\le k,\ell\le m.
\tag{A.1}
$$
We want to obtain conditions under which we can find real analytic solutions
$u$ to
$$
X_k u=0,\quad 1\le k\le m,
\tag{A.2}
$$
on a neighborhood of a given point $p\in\Cal{O}$.  We proceed as follows.
Say
$$
X_k=\sum\limits_j a_{kj}(x)\, \frac{\pa}{\pa x_j}.
\tag{A.3}
$$
On a neighborhood $\Omega$ of $p$ in $\CC^n$ set
$$
Z_k=\sum\limits_j a_{kj}(z)\, \frac{\pa}{\pa z_j},
\tag{A.4}
$$
with $a_{kj}(z)$ holomorphic extensions of $a_{kj}(x)$.  Solving (A.2) is
equivalent to finding a holomorphic solution $u$ to
$$
Z_k u=0,\quad 1\le k\le m,
\tag{A.5}
$$
on a neighborhood of $p$ in $\CC^n$.  Note that (A.1) implies
$$
[Z_k,Z_\ell]=0,\quad 1\le k,\ell\le m.
\tag{A.6}
$$

Our next step involves passing to real vector fields on $\Omega\subset\CC^n
\approx\RR^{2n}$.  Generally, if 
$$
Z=\sum\limits_j a_j(z)\, \frac{\pa}{\pa z_j},
\tag{A.7}
$$
set
$$
a_j(z)=f_j(z)+ig_j(z),
\tag{A.8}
$$
with $f_j$ and $g_j$ real valued, and then set
$$
\Phi(Z)=Y=\sum\limits_j \Bigl(f_j \frac{\pa}{\pa x_j}+g_j \frac{\pa}{\pa y_j}
\Bigr).
\tag{A.9}
$$
If $Z$ is a holomorphic vector field, i.e., if (A.7) holds with $a_j(z)$ 
holomorphic, we say $Y=\Phi Z$ is a {\it real-holomorphic vector field}.
Our first lemma holds whether or not the coefficients of $Z$ are holomorphic.

\proclaim{Lemma A.1}  If $a_j\in C(\Omega)$ in (A.6) and $Y=\Phi(Z)$, then
$$
u\ \text{ holomorphic}\Longrightarrow Zu=Yu.
\tag{A.10}
$$
\endproclaim

The proof is a straightforward calculation, making use of
$$
\frac{\pa u}{\pa z_j}=\frac{\pa u}{\pa x_j}=
\frac{1}{i} \frac{\pa u}{\pa y_j}.
\tag{A.11}
$$

The following is special to holomorphic vector fields, namely that $\Phi$
preserves the Lie bracket when applied to such vector fields.

\proclaim{Lemma A.2} If also $W=\sum b_j(z)\, \pa/\pa z_j$, then
$$
a_j,b_j\ \text{ holomorphic}\Longrightarrow \Phi[Z,W]=[\Phi Z,\Phi W].
\tag{A.12}
$$
\endproclaim

Again the proof is a straightforward (though slightly tedious) calculation.

It follows that if $X_k$ and $Z_k$ are as in (A.3)--(A.4), and if
$$
Y_k=\Phi Z_k=\sum\limits_j \Bigl( f_{kj} \frac{\pa}{\pa x_j}
+g_{kj} \frac{\pa}{\pa y_j}\Bigr),\quad f_{kj}=\Re a_{kj},\ g_{kj}=\Im a_{kj},
\tag{A.13}
$$
then
$$
[X_k,X_\ell]=0\Rightarrow [Z_k,Z_\ell]=0\Rightarrow [Y_k,Y_\ell]=0,\quad
1\le k,\ell\le m.
\tag{A.14}
$$

The complex structure on $\CC^n$ produces a complex structure on the space
of real vector fields on $\Omega$, defined by
$$
J\, \frac{\pa}{\pa x_j}=\frac{\pa}{\pa y_j},\quad 
J\, \frac{\pa}{\pa y_j}=-\frac{\pa}{\pa x_j}.
\tag{A.15}
$$
Note that if $Z$ has the form (A.7), then
$$
\Phi(iZ)=J\, \Phi(Z).
\tag{A.16}
$$
In particular, if $Y_k$ are as in (A.13),
$$
[Y_k,JY_\ell]=0=[JY_k,JY_\ell],\quad 1\le k,\ell\le m.
\tag{A.17}
$$

One advantage of using the real vector fields $Y_k$ on $\Omega$ is that they
generate local flows $\Cal{F}^t_{Y_k}$ on $\Omega$.  In this context,
the following results are very useful.

Suppose $Y$ is a real-holomorphic vector field on $\Omega$.  It follows from
(A.16) that so is $JY$, and $Y$ and $JY$ commute.  Thus so do the local
flows $\Cal{F}^s_Y$ and $\Cal{F}^t_{JY}$.  The following gives important
information on how these flows fit together.

\proclaim{Proposition A.3}  If $Y$ is a real-holomorphic vector field on
$\Omega$, then, for each $z\in\Omega$,
$$
\Cal{F}^s_Y \Cal{F}^t_{JY}(z)\ \text{ is holomorphic in }\ s+it.
\tag{A.18}
$$
\endproclaim
\demo{Proof}  Denote the 2-parameter orbit in  (A.18) by $\varphi(s,t)$.
By commutativity we also have
$$
\varphi(s,t)=\Cal{F}^t_{JY} \Cal{F}^s_Y(z).
\tag{A.19}
$$
It follows that
$$
\frac{\pa\varphi}{\pa s}=Y(\varphi(s,t)),\quad
\frac{\pa\varphi}{\pa t}=JY(\varphi(s,t)),
\tag{A.20}
$$
and hence $\pa \varphi/\pa t=J\, \pa\varphi/\pa s$, which gives the 
asserted holomorphicity.
\enddemo

The following is an important complement.

\proclaim{Proposition A.4} If $Y$ is a real-holomorphic vector field on 
$\Omega$, then $\Cal{F}^t_Y$ is a local group of holomorphic maps.
\endproclaim
\demo{Proof} The claim is equivalent to the assertion that
$$
\Cal{F}^t_{Y\#}\circ J=J\circ \Cal{F}^t_{Y\#},
\tag{A.21}
$$
where, given a diffeomorphism $F$, $F_{\#}$ is the induced operator on 
vector fields.  One has
$$
\frac{d}{dt} \Cal{F}^t_{Y\#}\, W\Bigr|_{t=0}=[Y,W];
\tag{A.22}
$$
cf.~(8.3) in Chapter I of [T].  Hence
$$
\frac{d}{dt} \Bigl( \Cal{F}^t_{Y\#}\circ J-J\circ\Cal{F}^t_{Y\#}\Bigr)W
\Bigr|_{t=0}=[Y,JW]-J[Y,W].
\tag{A.23}
$$
If $Y=\Phi Z$ with $Z$ a holomorphic vector field, as in (A.7)--(A.9), then
a calculation using
$$
\frac{\pa f_j}{\pa x_\ell}=\frac{\pa g_j}{\pa y_\ell},\quad
\frac{\pa f_j}{\pa y_\ell}=-\frac{\pa g_j}{\pa x_\ell},
\tag{A.24}
$$
shows that, for any vector field $W$,
$$
[Y,JW]-J[Y,W]=0,
\tag{A.25}
$$
so the quantity (A.23) vanishes.  More generally,
$$
\aligned
\frac{d}{dt}\bigl(\Cal{F}^t_{Y\#}\circ J-J\circ \Cal{F}^t_{Y\#}\bigr)W
&=\Cal{F}^t_{Y\#}[Y,JW]-J\Cal{F}^t_{Y\#}[Y,W] \\
&=(\Cal{F}^t_{Y\#}\circ J-J\circ \Cal{F}^t_{Y\#})[Y,W],
\endaligned
\tag{A.26}
$$
the latter identity by (A.25).  An iteration gives
$$
\Bigl(\frac{d}{dt}\Bigr)^\ell (\Cal{F}^t_{Y\#}\circ J-J\circ \Cal{F}^t_{Y\#})
W=(\Cal{F}^t_{Y\#}\circ J-J\circ \Cal{F}^t_{Y\#})(\Cal{L}^\ell_Y W),
\tag{A.27}
$$
where $\Cal{L}_YW=[Y,W]$.  In particular, for all $\ell\in\ZZ^+$,
$$
\Bigl(\frac{d}{dt}\Bigr)^\ell (\Cal{F}^t_{Y\#}\circ J-J\circ\Cal{F}^t_{Y\#})
W\bigr|_{t=0}=0.
\tag{A.28}
$$
In the current context, $\Cal{F}^t_Y$ and all its derived quantities are
real analytic in $t$ (as a consequence of Proposition A.3), so (A.21)
follows from (A.28).
\enddemo

We proceed to find solutions to (A.2), under appropriate hypotheses.
For notational simplicity, assume $p$ is the origin; $p=0\in\RR^n\subset
\CC^n$.  Suppose
$$
V\ \text{ is a linear subspace of }\ \RR^n,\ \text{ of dimension }\ n-m,
\tag{A.29}
$$
and let
$$
\widetilde{V}\ \text{ be the complexification of }\ V,
\tag{A.30}
$$
so $\widetilde{V}$ is a complex subspace of $\CC^n$, of complex dimension 
$n-m$ (hence real dimension $2n-2m$).  Let $v$ be a real analytic function
on a neighborhood $\Cal{U}$ of $0$ in $V$, extended to a holomorphic
function on a neighborhood $\widetilde{\Cal{U}}$ of $0$ in $\widetilde{V}$.
We assume
$$
\{Y_k,JY_k:1\le k\le m\}\ \text{ is transverse to }\ \widetilde{V},
\tag{A.31}
$$
on $\widetilde{\Cal{U}}$.  In particular,
$$
\CC^n=\RR\text{-linear span of }\ \widetilde{V}\ \text{ and }\
\{Y_k(0),JY_k(0):1\le k\le m\}.
\tag{A.32}
$$
Conversely, if (A.32) holds, then (A.31) holds, possibly with 
$\widetilde{\Cal{U}}$ shrunken.  In such a case, we can set
$$
u(\zeta)=v(z),\ \text{ for }\ z\in\widetilde{\Cal{U}},\ \ 
\zeta=\Cal{F}^{s_1}_{Y_1}\Cal{F}^{t_1}_{JY_1}\cdots
\Cal{F}^{s_m}_{Y_m}\Cal{F}^{t_m}_{JY_m}(z),
\tag{A.33}
$$
and see that $u$ is holomorphic on a neighborhood of $0$ in $\CC^n$ 
and solves
$$
Y_ku=JY_ku=0,\quad 1\le k\le m.
\tag{A.34}
$$
Hence, by Lemma A.1, (A.5) holds, hence, possibly shrinking $\Cal{U}$,
we have
$$
X_ku=0,\quad 1\le k\le m,\quad u\bigr|_{\Cal{U}}=v.
\tag{A.35}
$$

A classic example to which this construction applies arises in the real
analytic case of the Newlander-Nirenberg theorem.  In this setting, one
has $n=2m$ and takes $\xi_j=x_j+ix_{j+m},\ 1\le j\le m$, and
$$
X_k=\frac{\pa}{\pa\overline{\xi}_k}-\sum\limits_{\ell=1}^m b_{k\ell}(x)\, 
\frac{\pa}{\pa \xi_\ell},\quad 1\le k\le m,\quad b_{k\ell}(0)=0.
\tag{A.36}
$$
These vector fields arise from an almost complex structure $J_0$ on $\Cal{O}
\subset\RR^n$, and the integrability condition is that they commute, i.e., 
that (A.1) holds.  Then a function $u$ on $\Cal{O}$ is holomorphic with 
respect to this almost complex structure if and only if (A.2) holds, and
the theorem is that if (A.1) holds then there are $m$ such functions forming
a local coordinate system, in a neighborhood of $0$.  In this case we have
$$
X_k(0)=\frac{\pa}{\pa\overline{\xi}_k}=\frac{1}{2}\frac{\pa}{\pa x_k}
+\frac{i}{2}\frac{\pa}{\pa x_{m+k}},\quad 1\le k\le m,
\tag{A.37}
$$
hence
$$
Z_k(0)=\frac{1}{2}\frac{\pa}{\pa z_k}+\frac{i}{2} \frac{\pa}{\pa z_{m+k}},
\tag{A.38}
$$
so
$$
Y_k(0)=\frac{1}{2}\frac{\pa}{\pa x_k}+\frac{1}{2} \frac{\pa}{\pa y_{m+k}},
\tag{A.39}
$$
and
$$
JY_k(0)=\frac{1}{2}\frac{\pa}{\pa y_k}-\frac{1}{2}\frac{\pa}{\pa x_{m+k}}.
\tag{A.40}
$$

Let us take for $V\subset\RR^n$ the space
$$
V=\{x\in\RR^n:x_{m+1}=\cdots=x_{2m}=0\},
\tag{A.41}
$$
so
$$
\widetilde{V}=\{x+iy\in\CC^n:x_{m+1}=\cdots=x_{2m}=y_{m+1}=\cdots=y_{2m}=0\},
\tag{A.42}
$$
which is spanned over $\RR$ by
$$
\Bigl\{ \frac{\pa}{\pa x_j},\ \frac{\pa}{\pa y_j}:1\le j\le m\Bigr\}.
\tag{A.43}
$$
It is clear that if $Y_k(0)$ and $JY_k(0)$ are given by (A.39)--(A.40), then
(A.32) holds, so we have solutions to (A.35) in this case, for some 
neighborhood $\Cal{U}$ of $0$ in $V$, and arbitrary real analytic $v$ on
$\Cal{U}$.  This provides enough $J_0$-holomorphic functions on a 
neighborhood of $0$ in $\RR^n$ to yield a coordinate system.  In this fashion 
the real analytic case of the Newlander-Nirenberg theorem is proven.

$$\text{}$$
{\bf B. The case of two-dimensional leaves}
\newline {}\newline

Here we put ourselves in the setting of \S{3}, and take the Lipschitz bundle
$\Cal{E}$ to have fiber dimension $k=2m=2$.  We assume $\Cal{E}$ has a 
complex structure $J$, pulled back as in \S{3} to a complex structure 
$J^0\in\End(\Cal{E}_0)$, where $\Cal{E}_0\subset TU_0$ is the bundle 
spanned by $\pa/\pa t_1,\ \pa/\pa t_2$.  Here $U_0\subset\RR^n$ is an open
set with coordinates $(t,z),\ t\in\RR^2,\ z\in\RR^{n-2}$.  We assume 
$$
J\in C^r(U_1),
\tag{B.1}
$$
with $r\in(0,1)$, in which case
$$
J^0\in C^r(U_0).
\tag{B.2}
$$
We can represent $J^0=J^0(t,z)$ as a $2\times 2$ matrix valued function
of $(t,z)$.  Making a preliminary change of coordinates 
$$
(t,z)\mapsto
(A(z)t,z),
\tag{B.3}
$$
where $A(z)$ is a $G\ell(2,\RR)$-valued function of the
same type of regularity as $J^0$ in (B.2), we can arrange that
$$
J^0(0,z)=\pmatrix 0 & -1\\ 1 & 0\endpmatrix ,
\tag{B.4}
$$
for all $z$.

In order to implement the classical method of finding isothermal coordinates,
we impose a family of Riemannian metric tensors on $t$-space, depending on
$z$ as a parameter, $(g_{ij}(t,z)),\ 1\le i,j\le 2$.  Arrange that $J^0(t,z)$
is an isometry on $T_t \RR^2$ with respect to the induced inner product,
for each $(t,z)$.  One could, for example, start with the standard flat
metric $(\delta_{ij})$ and average with respect to the $\ZZ/(4)$-action
generated by $J^0$.  We then obtain
$$
g_{ij}\in C^r(U_0),
\tag{B.5}
$$
when (B.2) holds, and we can arrange that
$$
g_{ij}(0,z)=\delta_{ij}.
\tag{B.6}
$$

Let $D=\{t\in\RR^2:t_1^2+t_2^2<1\}$.  We want to find a harmonic function
$u_1$ on $D$ equal to $t_1$ on $\pa D$ (and depending on the parameter $z$).
Thus, with $a^{ij}(t,z)=g(t,z)^{1/2} g^{ij}(t,z)$, where $(g^{ij})$ is the
inverse of the matrix $(g_{ij})$ and $g$ its determinant, we want to solve
$$
\pa_i a^{ij}(t,z)\pa_ju_1=0\ \text{ on }\ D,\quad u_1\bigr|_{\pa D}=t_1,
\tag{B.7}
$$
where $\pa_i=\pa/\pa t_i,\ i=1,2$.  Without changing notation, we dilate the 
$t$-coordinates, and we can assume
$$
a^{ij}(0,z)=\delta^{ij},\quad \|a^{ij}(\cdot,z)-\delta^{ij}\|_{C^r(\Dbar)}
\le\eta,
\tag{B.8}
$$
where $\eta>0$ is a sufficiently small quantity.  Let us write (B.7) as
$$
(\Delta+R_z)u_1=0,\quad u_1\bigr|_{\pa D}=t_1,
\tag{B.9}
$$
where
$$
R_zu_1=\pa_i r^{ij}(t,z)\pa_j u_1,\quad r^{ij}(t,z)=a^{ij}(t,z)-\delta^{ij},
\tag{B.10}
$$
and
$$
\Delta=\frac{\pa^2}{\pa t_1^2}+\frac{\pa^2}{\pa t_2^2}.
$$

To establish solvability of (B.9), when $\eta$ in (B.8) is small enough, 
note that it is equivalent to the following equation for $v=u_1-t_1$:
$$
(\Delta+R_z)v=-R_zt_1,\quad v\bigr|_{\pa D}=0,
\tag{B.11}
$$
hence to the equation
$$
(I+GR_z)v=-GR_z t_1,
\tag{B.12}
$$
where $G$ is the solution operator to the Poisson problem for $\Delta$ on 
$D$, with the Dirichlet boundary condition.  Such $G$ has the property
$$
G:C^{r-1}_*(\Dbar)\longrightarrow C^{r+1}(\Dbar),\quad 0<r<1;
\tag{B.13}
$$
cf.~[T], Chapter 13, (8.54)--(8.55).  Hence
$$
\aligned
\|GR_zf\|_{C^{r+1}(\Dbar)}&\le C\|r^{ij}(\cdot,z)\pa_jf\|_{C^r(\Dbar)} \\
&\le C\|r^{ij}(\cdot,z)\|_{C^r(\Dbar)} \|f\|_{C^{r+1}(\Dbar)},
\endaligned
\tag{B.14}
$$
so if $\eta$ is small enough, the operator norm of $GR_z$ on $C^{r+1}(\Dbar)$
is $\le 1/2$, so $I+GR_z$ in (B.12) is invertible on $C^{r+1}(\Dbar)$, and
we have a unique solution $v$, satisfying
$$
\aligned
\|v\|_{C^{r+1}(\Dbar)}&\le C\|GR_zt_1\|_{C^{r+1}(\Dbar)} \\
&\le C\|r^{i1}(\cdot,z)\|_{C^r(\Dbar)} \\ &\le C\eta.
\endaligned
\tag{B.15}
$$

We now have $u_1=t_1+v$.  The standard construction of the harmonic conjugate
$u_2$, satisfying
$$
du_2=(J^0)^t du_1,\quad u_2(0,z)=0,
\tag{B.16}
$$
gives
$$
\|u_2-t_2\|_{C^{r+1}(\Dbar)}\le C\eta,
\tag{B.17}                  
$$
and taking $u=u_1+iu_2$, we have a local holomorphic coordinate system on
each leaf $z=z_0$, if $\eta$ is small enough.

We now want to determine how smooth $u_i(t,z)$ are in $z$, first in the
case $i=1$.  So pick points $z$ and $z'$ and set $w=u_1(t,z)-u_1(t,z')$.
Hence
$$
\Delta w=-R_{z}u_1(\cdot,z)+R_{z'}u_1(\cdot,z'),\quad 
w\bigr|_{\pa D}=0,
\tag{B.18}
$$
or alternatively
$$
(\Delta+R_{z})w=-(R_{z}-R_{z'})u_1(\cdot,z'),\quad 
w\bigr|_{\pa D}=0.
\tag{B.19}
$$
An argument similar to (B.11)--(B.14) yields, for $s\in (0,r]$,
$$
\|w\|_{C^{1+s}(\Dbar)}\le C\|r^{ij}(\cdot,z)-r^{ij}(\cdot,z')
\|_{C^s(\Dbar)} \|u_1(\cdot,z')\|_{C^{s+1}(\Dbar)}.
\tag{B.20}
$$
We already have a bound on $u_1=t_1+v$ from (B.15).  As for the other factor
on the right side of (B.20), we can use the elementary estimate
$$
\|f\|_{C^s(\Dbar)}\le C \|f\|_{C^r(\Dbar)}^{s/r}
\|f\|_{L^\infty(\Dbar)}^{1-s/r},
\tag{B.21}
$$
valid for $s\in [0,r]$, to deduce that
$$
\|u_1(\cdot,z)-u_1(\cdot,z')\|_{C^{1+s}(\Dbar)}
\le C_s |z-z'|^{r-s},\quad 0<s\le r,
\tag{B.22}
$$
given the latter alternative in hypothesis (B.1).  The construction of
$u_2$ via (B.16) then yields
$$
\|u_2(\cdot,z)-u_2(\cdot,z')\|_{C^{1+s}(\Dbar)}
\le C_s |z-z'|^{r-s},\quad 0<s\le r.
\tag{B.23}
$$

Thus if we set
$$
u_{ij}=\frac{\pa u_i}{\pa t_j},\quad 1\le i,j\le 2,
\tag{B.24}
$$
we have
$$
\|u_{ij}(\cdot,z)-u_{ij}(\cdot,z')\|_{C^s(\Dbar)}\le C_s |z-z'|^{r-s},
\quad 0<s\le r,
\tag{B.25}
$$
and hence, taking respectively $s=r$ and $s=\delta$, close to $0$, we have
$$
\aligned
|u_{ij}(t,z)-u_{ij}(t',z')|&\le |u_{ij}(t,z)-u_{ij}(t',z)|
+|u_{ij}(t',z)-u_{ij}(t',z')| \\ &\le C|t-t'|^r+C_\delta |z-z'|^{r-\delta}.
\endaligned
\tag{B.26}
$$
If we reverse the coordinate transformation (B.3), this estimate remains 
valid.

We obtain a CR function on an open set in $U_1$ by composing $u=u_1+iu_2$
with the inverse of the bi-Lipschitz map $G$, given in (3.1):
$$
\tilde{u}=u\circ G^{-1}.
\tag{B.27}
$$
We have
$$
Y_j\tilde{u}=\frac{\pa u}{\pa t_j}\circ G^{-1},\quad j=1,2,
\tag{B.28}
$$
where $\{Y_1,Y_2\}$ is the Lipschitz frame field for $\Cal{E}$ that pulls back
to $\{\pa/\pa t_1,\pa/\pa t_2\}$.  It follows that
$$
\tilde{u},\ Y_j\tilde{u}\in C^{r-\delta},\quad \forall\, \delta>0.
\tag{B.29}
$$

While one cannot take $\delta=0$ in (B.26), one can improve the estimate,
as follows.  First, using
$$
G:H^{-1,p}(D)\longrightarrow H^{1,p}(D),\quad 1<p<\infty,
\tag{B.30}
$$
an argument parallel to (B.11)--(B.14) gives, in place of (B.20),
$$
\aligned
\|w\|_{H^{1,p}(D)}&\le C\|(R_z-R_{z'})u_1(\cdot,z')\|_{H^{-1,p}(D)} \\
&\le C\|r^{ij}(\cdot,z)-r^{ij}(\cdot,z')\|_{L^\infty(D)}
\|u_1(\cdot,z')\|_{H^{1,p}(D)}.
\endaligned
\tag{B.31}
$$
Then one can exploit the following local regularity result.  Suppose 
$\omega(h)$ is a modulus of continuity satisfying the Dini condition:
$$
\int_0^{1/2} \frac{\omega(h)}{h}\, dh<\infty.
\tag{B.32}
$$
Then, with $D_{1/2}=\{t:t_1^2+t_2^2<1/4\}$,
$$
\aligned
\|w\|_{C^1(D_{1/2})}
&\le C\|w\|_{H^{1,p}(D)}+C\|(R_z-R_{z'})u_1(\cdot,z')\|_{C^{-1,\omega}(D)} \\
&\le C\|w\|_{H^{1,p}(D)}+C\|(r^{ij}(\cdot,z)-r^{ij}(\cdot,z'))\pa_i
u_1(\cdot,z')\|_{C^\omega(\Dbar)} \\
&\le C\|w\|_{H^{1,p}(D)}+C\|r^{ij}(\cdot,z)-r^{ij}(\cdot,z')
\|_{C^\omega(\Dbar)} \|u_1(\cdot,z')\|_{C^{r+1}(\Dbar)}.
\endaligned
\tag{B.33}
$$
In view of (B.31), and the previous estimates on $u_1(\cdot,z')$, we have
$$
\|u_1(\cdot,z)-u_1(\cdot,z')\|_{C^1(D_{1/2})}\le
C\|r^{ij}(\cdot,z)-r^{ij}(\cdot,z')\|_{C^\omega(\Dbar)}.
\tag{B.34}
$$

To estimate the right side of (B.34), we replace (B.21) by the following.
Suppose
$$
|f(x-y)|\le C|x-y|^r,\quad |f(x)-f(y)|\le C\delta.
\tag{B.35}
$$
Then
$$
|f(x)-f(y)|\le C\sigma_r(\delta)\, \omega(|x-y|),
\tag{B.36}
$$
where
$$
\sigma_r(\delta)=\sup\limits_{h\in (0,1]}\, 
\frac{\min(\delta,h^r)}{\omega(h)}.
\tag{B.37}
$$
Of course, we pick $\omega(h)$ decreasing to $0$ as $h\searrow 0$ more slowly
than $h^r$ for any $r>0$, for example,
$$
\omega(h)=\Bigl(\log \frac{1}{h}\Bigr)^{-1-a},
\tag{B.38}
$$
for some $a>0$, so that $h^r/\omega(h)$ is $\nearrow$ on $h\in (0,1/2]$.  
In such a case,
$$
\sigma_r(\delta)\approx \frac{\delta}{\omega(\delta^{1/r})}.
\tag{B.39}
$$

We deduce that, under the hypothesis (B.1),
$$
\aligned
\|u_1(\cdot,z)-u_1(\cdot,z')\|_{C^1(D_{1/2})}&\le C\, \sigma_r(|z-z'|^r) \\
&\le C \frac{|z-z'|^r}{\omega(|z-z'|)}.
\endaligned
\tag{B.40}
$$
Hence we can supplement (B.25) with
$$
\|u_{ij}(\cdot,z)-u_{ij}(\cdot,z')\|_{C^0(D_{1/2})}
\le C \frac{|z-z'|^r}{\omega(|z-z'|)},
$$
and improve (B.26) to
$$
|u_{ij}(t,z)-u_{ij}(t',z')|\le C|t-t'|^r
+C \frac{|z-z'|^r}{\omega(|z-z'|)}.
$$
This in turn leads to an improvement in the modulus of continuity in (B.29), 
to
$$
|\tilde{u}(x)-\tilde{u}(x')|,\ |Y_j\tilde{u}(x)-Y_j\tilde{u}(x')|
\le C \frac{|x-x'|^r}{\omega(|x-x'|)}.
\tag{B.41}
$$

Next, we want to replace hypothesis (B.1) by
$$
J\in\Lip(U_1).
\tag{B.42}
$$
Then we replace $C^r$ by Lip in (B.2), (B.5), and (B.8), and we supplement 
(B.13) by
$$
G:C^0_*(\Dbar)\longrightarrow C^2_*(\Dbar).
\tag{B.43}
$$
Thus (B.14) is modified to
$$
\|GR_zf\|_{C^2_*(\Dbar)}\le C\|r^{ij}(\cdot,z)\|_{C^1_*(\Dbar)}
\|f\|_{C^2_*(\Dbar)},
\tag{B.44}
$$
which leads to the existence of isothermal coordinates $u_1, u_2$, satisfying
$$
\|u_i-t_i\|_{C^2_*(\Dbar)}\le C\eta.
\tag{B.45}
$$
Analogues of (B.18)--(B.23) hold.  We need to replace (B.21) by the 
interpolation inequality
$$
\|f\|_{C^s_*(\Dbar)}\le C \|f\|^s_{C^1_*(\Dbar)}\|f\|^{1-s}_{C^0_*(\Dbar)},
\tag{B.46}
$$
valid for $s\in [0,1]$, and then we get
$$
\|u_i(\cdot,z)-u_i(\cdot,z')\|_{C^{1+s}_*(\Dbar)}
\le C_s |z-z'|^{1-s},\quad 0<s\le 1.
\tag{B.47}
$$
Keep in mind that $C^{1+s}_*(\Dbar)=C^{1+s}(\Dbar)$ for $0<s<1$.  Similarly,
in place of (B.25), we have
$$
\|u_{ij}(\cdot,z)-u_{ij}(\cdot,z')\|_{C^s_*(\Dbar)}\le C_s|z-z'|^{1-s},
\quad 0<s\le 1.
\tag{B.48}
$$

Consequently (B.26) is modified as follows.  First, since elements of
$C^1_*(\Dbar)$ have a log-Lipschitz modulus of continuity, we have
$$
|u_{ij}(t,z)-u_{ij}(t',z)|\le C|t-t'|\, \log\frac{1}{|t-t'|},
\tag{B.49}
$$
for $t,t'\in\Dbar$.  On the other hand, (B.34) still applies, and we can
take $r=1$ in (B.37) to obtain
$$
\|u_{ij}(\cdot,z)-u_{ij}(\cdot,z')\|_{C^0(D_{1/2})}\le
C\sigma_1(|z-z'|),
\tag{B.50}
$$
where
$$
\sigma_1(\delta)=\sup\limits_{h\in(0,1]}\, \frac{\min(\delta,h)}{\omega(h)}.
\tag{B.51}
$$
We thus obtain, in place of (B.40), the modulus of continuity estimate
$$
|u_{ij}(t,z)-u_{ij}(t',z')|\le C|t-t'|\, \log\frac{1}{|t-t'|}
+C\sigma_1(|z-z'|),
\tag{B.52}
$$
for $t,t'\in D_{1/2}$.  This in turn leads to
$$
|\tilde{u}(x)-\tilde{u}(x')|,\ |Y_j\tilde{u}(x)-Y_j\tilde{u}(x')|
\le C\sigma^\#(|x-x'|),
\tag{B.53}
$$
where
$$
\sigma^\#(\delta)=\max\Bigl(\sigma_1(\delta),\delta\, \log\frac{1}{\delta}
\Bigr),\quad \text{for }\ 0<\delta\le 1.
\tag{B.54}
$$
If $\omega(h)$ is given by (B.38), we have
$$
\sigma^\#(\delta)=\delta\, \Bigl(\log \frac{1}{\delta}\Bigr)^{1+a}.
\tag{B.55}
$$

We formally state the main conclusion of this appendix.  Since the result
is local, we may as well take $\Omega$ to be an open set in some Euclidean
space.

\proclaim{Proposition B.1} Let $\Omega$ have a Lipschitz, Levi-flat 
CR-structure, with leaves tangent to $\Cal{E}$ of real dimension two.  Then 
each $p\in\Omega$ has a neighborhood $U$ on which there is a CR-function
$$
\tilde{u}:U\longrightarrow \CC,
\tag{B.56}
$$
which is a holomorphic diffeomorphism on each leaf, intersected with $U$,
into $\CC$, with the following regularity.  For any $a>0$, and any
Lipschitz section $Y$ of $\Cal{E}$,
$$
|\tilde{u}(x)-\tilde{u}(x')|,\ |Y\tilde{u}(x)-Y\tilde{u}(x')|\le
C_a\, |x-x'|\, \Bigl(\log \frac{1}{|x-x'|}\Bigr)^{1+a},
\tag{B.57}
$$
for $x,x'\in U,\ |x-x'|\le 1/2$.
\endproclaim

$\text{}$ \newline
{\smc Remark}. Since a tool in the analysis of the Lipschitz CR-structures
was an analysis of families of much less regular almost complex structures,
it is worth mentioning the fundamental work of Ahlfors and Bers [AB] on
the endpoint case, involving merely $L^\infty$ almost complex structures.
See also [A] and [D] for treatments; the latter article also discusses
dependence on parameters.  In such a case the $C^1$ regularity collapses to
H{\"o}lder continuity, and it does not seem that techniques used there lead 
to an improvement of Proposition B.1.

$$\text{}$$
{\bf References}

$\text{}$
\roster
\item"[A]" L.~Ahlfors, Lectures on Quasiconformal Mappings, Wadsworth, 1987.
\item"[AB]" L.~Ahlfors and L.~Bers, Riemann's mapping theorem for variable 
metrics, Annals of Math. 72 (1960), 385--404.
\item"[AH]" A.~Andreotti and C.D.~Hill, Complex characteristic coordinates
and the tangential Cauchy-Riemann equations, Ann. Scuola Norm. Sup. Pisa
26 (1972), 299--324.
\item"[Bog]" A.~Boggess, CR Manifolds and the Tangential Cauchy-Riemann
Complex, CRC Press, Boca Raton, Florida, 1991.
\item"[D]" A.~Douady, Le th{\'e}or{\`e}me d'integrabilit{\'e} des 
structures presque complexes, pp.~307--324 in ``The Mandelbrot Set, Theme and 
Variations,'' Tan Lei (ed.), Cambridge Univ. Press, 2000.
\item"[Fr]" M.~Freeman, The Levi form and local complex foliations,
Proc. AMS 57 (1976), 368--370.
\item"[Ha]" P.~Hartman, Frobenius theorem under Carath{\'e}odory type
conditions, J. Diff. Eqns. 7 (1970), 307--333.
\item"[HT]" C.D.~Hill and M.~Taylor, Integrability of rough almost complex
structures, J. Geom. Anal. 13 (2003), 163--172.
\item"[Ho]" L.~H{\"o}rmander, The Frobenius-Nirenberg theorem, Arkiv f{\"o}r
Matematik 5 (1964), 425--432.
\item"[LM]" C.~Lebrun and L.~Mason, Zoll manifolds and complex surfaces,
J. Diff. Geom. 61 (2002), 453--535.
\item"[LC]" T.~Levi-Civita, Sulle funzione di due o pi{\`u} variabli
complesse, Rend. Acc. Lincei 14 (1905), 492--499.
\item"[M]" B.~Malgrange, Sur l'int{\'e}grabilit{\'e} des structures
presque-complexes, Symposia Math., Vol. II (INDAM, Rome, 1968), Academic 
Press, London, 289--296, 1969.
\item"[NN]" A.~Newlander and L.~Nirenberg, Complex coordinates in almost
complex manifolds, Ann. of Math. 65 (1957), 391--404.
\item"[NW]" A.~Nijenhuis and W.~Woolf, Some integration problems in
almost-complex manifolds, Ann. of Math. 77 (1963), 424--489.
\item"[Ni]" L.~Nirenberg, A complex Frobenius theorem, Seminars on Analytic
Functions I, 172--189.  Institute for Advanced Study, Princeton, 1957.
\item"[Pin]" S.~Pinchuk, CR-transformations of real manifolds in $\CC^n$,
Indiana Univ. Math. J. 41 (1992), 1--15.
\item"[Som]" F.~Sommer, Komplex-analytische Blaetterung reeler
Manifaltigkeiten im $\CC^n$, Math. Annalen 136 (1958), 111-113.
\item"[T]" M.~Taylor, Partial Differential Equations, Vols.~1--3, 
Springer-Verlag, New York, 1996.
\endroster

$$\text{}$$

{\smc Department of Mathematics, Stony Brook University, Stony Brook,
New York 11794}

{\it E-mail address:} {\tt dhill\@math.sunysb.edu}

$\text{}$

{\smc Department of Mathematics, University of North Carolina, Chapel Hill,
North Carolina 27599}

{\it E-mail address:} {\tt met\@math.unc.edu}

\end